\documentclass[a4paper, 12pt]{article}
\usepackage{amsfonts}
\usepackage{amssymb}
  \usepackage{color}

\input epsf

\newif\ifscrf\scrftrue
\ifx\footscrfont\nullfont
  \scrffalse
\fi
\newif\iffn\fnfalse

\ifscrf
  \font\footscrfont=rsfs10
  \font\footbbbfont=msbm10
  
\fi

\ifscrf
  
\else
  
\fi

\let\useblackboard=\iftrue
\newfam\black
\useblackboard
\font\blackboard=msbm10 scaled \magstep1
\font\blackboards=msbm7
\font\blackboardss=msbm5
\textfont\black=\blackboard
\scriptfont\black=\blackboards
\scriptscriptfont\black=\blackboardss

\else

\fi

\def\yboxit#1#2{\vbox{\hrule height #1 \hbox{\vrule width #1
\vbox{#2}\vrule width #1 }\hrule height #1 }}
\def\fillbox#1{\hbox to #1{\vbox to #1{\vfil}\hfil}}
\def\ybox{{\lower 1.3pt \yboxit{0.4pt}{\fillbox{8pt}}\hskip-0.2pt}}
\def\CC{{\cal C}}

\def\inbar{\,\vrule height1.5ex width.4pt depth0pt}
\font\cmss=cmss10 \font\cmsss=cmss10 at 7pt

\def\IZ{\relax\ifmmode\mathchoice
{\hbox{\cmss Z\kern-.4em Z}}{\hbox{\cmss Z\kern-.4em Z}}
{\lower.9pt\hbox{\cmsss Z\kern-.4em Z}}
{\lower1.2pt\hbox{\cmsss Z\kern-.4em Z}}\else{\cmss Z\kern-.4em
Z}\fi}
\def\IB{\relax{\rm I\kern-.18em B}}
\def\IC{{\relax\hbox{$\inbar\kern-.3em{\rm C}$}}}
\def\ID{\relax{\rm I\kern-.18em D}}
\def\IE{\relax{\rm I\kern-.18em E}}
\def\IF{\relax{\rm I\kern-.18em F}}
\def\IG{\relax\hbox{$\inbar\kern-.3em{\rm G}$}}
\def\IGa{\relax\hbox{${\rm I}\kern-.18em\Gamma$}}
\def\IH{\relax{\rm I\kern-.18em H}}
\def\II{\relax{\rm I\kern-.18em I}}
\def\IK{\relax{\rm I\kern-.18em K}}
\def\IP{\relax{\rm I\kern-.18em P}}
\def\IR{\relax{\rm I\kern-.18em R}}

\newif\ifamsf\amsftrue
\ifx\footbbbfont\nullfont
  \amsffalse
\fi

\ifamsf
  
\fi

\ifamsf
  \def\Goth#1{\iffn
      \mathchoice{\hbox{\footgthfont #1}}{\hbox{\footgthfont #1}}
      {\hbox{\smallfootgthfont #1}}{\hbox{\tinyfootgthfont #1}}\else
      \mathchoice{\hbox{\gthfont #1}}{\hbox{\gthfont #1}}
      {\hbox{\smallgthfont #1}}{\hbox{\tinygthfont #1}}\fi}
\else
  
  \def\Goth{\cal}
\fi

\def\Hom{\operatorname{Hom}}

\def\tr{\operatorname{tr}}

\newtheorem{thm}{Theorem}[section]
\newtheorem{lemma}[thm]{Lemma}

\newtheorem{cor}[thm]{Corollary}
\newtheorem{prop}[thm]{Proposition}
\newtheorem{defn}[thm]{Definition}

\newenvironment{proof}{{\sc Proof.}}{\ q.e.d. \vskip .2in}

\newtheorem{introthm}{Theorem}
\newtheorem{appendthm}{Theorem}
\newtheorem{appendlemma}[appendthm]{Lemma}

\def\cal{\mathcal}
\def\em{\it}
\def\r{\rho_0}
\def\mera{m_{\alpha}}

\def\a{{\alpha}}

\def\b{{\beta}}
\def\C{{\cal C}}

\def\l{{\lambda}}
\def\L{{\cal L}}

\def\M{{\cal M}}
\def\P{{\cal P}}
   
   \def\R{{\cal R}}
\def\s{{\sigma}}

\def\t{{\theta}}

\def\CC{{\mathbb C}}

\def\HH{{\mathbb H}}

\def\RR{{\mathbb R}}

\def\ZZ{{\mathbb Z}}

\def\ch {{\cosh}}
\def\sh {{\sinh}}

\def\dd{ {\partial}}

\def\pl{ pl}

\def\teich{{\rm Teich}}

\def\tr{\mathop{\rm Tr}}
\def\p_#1{ \frac{\partial \hfil}{\partial {#1}}}
\def\u#1 {\underline #1  }
\def\ua{ \underline{\alpha}}

\def\sl2c{{\Goth s}{\Goth l}(2,\CC)}
\def\Ad{{\mathop {\rm Ad}}}
\def\Hom{{\rm{Hom}}}
\def\Ax{\mathop{\rm Ax}}

\begin{document}
\renewcommand{\thefootnote}{\fnsymbol{footnote}}    

\begin{center}   
{\LARGE Lengths are coordinates for convex structures}

\vspace{5mm}
{\large Young-Eun Choi\footnote{Current Address:  Department of
    Mathematics, University of   California, 
Davis, CA 95616, USA.
E-mail: choiye@math.ucdavis.edu} 
 and Caroline Series\footnote{E-mail: cms@maths.warwick.ac.uk}}
\vspace{3mm}

{Mathematics Institute\\ 
University of Warwick\\
Coventry CV4 7AL, UK
}

\end{center}

\bibliographystyle{Plain}

\begin{abstract} Suppose that $N$ is a geometrically finite
  orientable hyperbolic $3$-manifold.  Let $\P(N,\ua)$ be the space of
  all geometrically finite hyperbolic structures on $N$ whose convex
  core is bent along a set $\ua$ of simple closed curves. We prove
  that the map which associates to each structure in $\P(N,\ua)$ the
  lengths of the curves in the bending locus $\ua$ is one-to-one. If
  $\ua$ is maximal, the traces of the curves in $\ua$ are local
  parameters for the representation space $\R(N)$.\\

\noindent{\em Key words: convex structure, cone manifold, bending
    lamination}
\end{abstract}

\renewcommand{\thefootnote}{\arabic{footnote}}                                  \setcounter{footnote}{0}   

\section{Introduction}

This paper is about the parameterization of convex structures on
hyperbolic $3$-manifolds.  We show that the space of structures whose
bending locus is a fixed set of closed curves $\ua = \{ \a_1,\ldots,
\a_n\}$ is parameterized by the hyperbolic lengths $\{l_{\a_1},\ldots,
l_{\a_n}\}$ and moreover, that when $\ua$ is maximal, the complex
lengths (or traces) of the same set of curves are local holomorphic
parameters for the ambient deformation space.

Suppose that $N(G)=\HH^3/G$ is a hyperbolic $3$-manifold such that $G$
is geometrically finite with non-empty regular set and such that its
convex core $V=V(G)$ has finite but non-zero volume.  The boundary of
$V$ is always a pleated surface, with bending locus a geodesic
lamination $\pl(G)$. Assume that $V(G)$ has {\em no rank-$2$ cusps}.
All our main results hold without this assumption, but writing down
the proofs in full generality does not seem to warrant all the
additional comments and notation entailed.  If $V$ has rank-$1$ cusps,
compactify $V$ by removing a horoball neighborhood of each cusp. The
interior of the resulting manifold $\bar N(G)$ is homeomorphic to
$\HH^3/G$. The boundary of $\bar N(G)$ consists of closed surfaces
obtained from $\dd V$ by removing horospherical neighborhoods of
matched pairs of punctures and replacing them with annuli, see
Figure~\ref{fig:annuli}. Denote by $\ua_P(G)$ the collection of core
curves of these annuli.

Now let $\bar N$ be a compact orientable $3$-manifold whose interior
$N$ admits a complete hyperbolic structure. Assume that the boundary
$\dd \bar N$ is non-empty and that it contains no tori. Let
$\underline \a = \{\a_1, \ldots, \a_m \}$ be a collection of disjoint,
homotopically distinct simple closed curves on $\dd \bar N$. In this
paper we investigate hyperbolic $3$-manifolds $\HH^3/G$ such that
$\bar N(G)$ is homeomorphic to $\bar N$ and such that the bending
locus of $\dd V$, together with the set $\ua_P(G)$ of core curves of
the annuli described above, consists exactly of the curves in $\ua$.
Denote by $\cal P^+(N,\ua)$ the space of all such structures,
topologized as a subset of the representation space $\cal R(N)
=Hom(\pi_1(N),SL(2,\CC))/SL(2,\CC)$, where $Hom(\pi_1(N),SL(2,\CC))$
is the space of homomorphisms from $\pi_1(N)$ to $SL(2,\CC)$ and
$SL(2,\CC)$ acts by conjugation. More generally, define the {\em
  pleating variety} $\P(N,\ua)$ to be the space of structures for
which the bending locus of $\dd V$ and $\ua_P(G)$ are contained in,
but not necessarily equal to, $\ua$.  We refer to a structure in
$\P(N,\ua)$ as a {\em convex structure} on $(\bar N,\ua)$. We
emphasize that structures in $\P(N,\ua)$ must have convex cores with
finite non-zero volume, thus excluding the possibility that the group
$G$ is Fuchsian.  

For $\a_i \in \ua$, let $\theta_{i}$ be the exterior bending angle
along $\a_i$, measured so that $\theta_{i}=0$ when $\a_i$ is contained
in a totally geodesic part of $\dd V$, and set $\theta_{i} = \pi$ when
$\a_i \in \ua_P(G)$.  Then $\cal P^+(N,\ua)$ is the subset of
$\P(N,\ua)$ on which $\theta_{i} > 0$ for all $i$.

In~\cite{BonO},  Bonahon and Otal gave necessary and sufficient
conditions for the existence of a convex structure with a given set of
bending angles. They show:
\begin{thm}[Angle parameterization]
\label{thm:boglobalpar}
Let $\Theta: \P^+(N,\ua) \rightarrow \RR^n$ be the map which
associates to each structure $\sigma$ the bending angles
$(\theta_1(\sigma),\ldots, \theta_n(\s))$ of the curves in the bending
locus $\ua=\{\a_1,\ldots,\a_n\}$. Then $\Theta$ is a diffeomorphism
onto a convex subset of $(0,\pi]^n$.
\end{thm}
Moreover, the image is entirely  specified by the topology of $\bar N$
and the curve  system $\ua$. Reformulating their conditions
topologically, we show in Theorem~\ref{thm:nonemptypleating} that
$\cal P^+(N,\ua)$ is non-empty exactly when the curves $\ua$ form a
{\em doubly incompressible system} on $\dd \bar N$, see
Section~\ref{subsec:topconds}.   Thus their result shows that,
provided these topological conditions  are satisfied,  $\P^+(N,\ua)$
is a submanifold of $\cal R(N)$ of real dimension equal to the number
of curves in $\ua$, and that the bending angles uniquely determine the
hyperbolic structure on $N$.

In this paper we prove an analogous parameterization theorem for
$\P(N,\ua)$ in which we replace the {\em angles} along the bending
lines by their {\em lengths}.
\renewcommand{\theintrothm}{\Alph{introthm}}
\begin{introthm}[Length parameterization]
\label{thm:globalpar}
Let $L: \P(N,\ua) \rightarrow \RR^n$ be the map which associates to
each structure $\sigma$ the hyperbolic lengths $(l_{\a_1}(\sigma),
\ldots,l_{\a_n}(\s))$ of the curves in the bending locus
$\ua=\{\a_1,\ldots,\a_n\}$. Then $L$ is an injective local diffeomorphism.
\end{introthm}

This follows from our stronger result that {\em any} combination of
lengths and angles also works:
\begin{introthm}[Mixed parameterization]
\label{cor:mixedparams}
For any ordering of the curves $\ua$ and for any $k$, the map $\sigma
\mapsto (l_{\alpha_1}(\s), \ldots, l_{\alpha_k}(\s),
\theta_{\alpha_{k+1}}(\s), \ldots,\theta_{\alpha_{n}}(\s))$ is 
an injective local  diffeomorphism on $\P(N,\ua)$.
\end{introthm}
\noindent{\bf Remark}
With some further work, one can show that $L$ is actually a
diffeomorphism onto its image. We hope to explore this, together with
some applications of the parameterization, elsewhere.

\medskip

It is known that in the neighborhood of a geometrically finite
representation, $\R(N)$ is a smooth complex variety of dimension equal
to the number of curves $d$ in a maximal curve system on $\dd \bar N$
(see Theorem~\ref{thm:smoothreps}).  Theorem~\ref{thm:globalpar}
follows from the following result on local parameterization:
\begin{introthm}[Local parameterization]
\label{thm:main2}
Let $\sigma_0 \in \P(N,\ua)$, where $\ua= \{\a_1,\ldots,\a_d\}$ is a
maximal curve system on $\dd \bar N$. Then the map $\L:\R(N)
\rightarrow \CC^d$ which associates to a structure $\sigma \in \cal
R(N)$ the complex lengths $(\lambda_{\a_1}(\sigma), \ldots,
\lambda_{\a_d}(\sigma))$ of the curves $\a_1,\ldots,\a_d$ is a local
diffeomorphism in a neighborhood of $\sigma_0$. 
\end{introthm}
To be precise, if $\s_0(\a_i)$ is parabolic, then $\lambda_{\a_i}(\s)$
must be replaced with $\tr \s(\a_i)$ in the definition of $\L$.  This
point, together with a precise definition of the complex length
$\lambda_{\a_i}(\sigma)$, is discussed in detail in
Section~\ref{subsec:complexlength}. It is not hard to show that when
restricted to $\cal P(N,\ua)$, the map $\cal L: \R(N) \rightarrow
\CC^d$ is real-valued and coincides with $L$.  We remark that the map
$\L$ is not globally non-singular; in fact we showed in~\cite{S2} that
if $G$ is quasifuchsian, then $\mbox{d}\L$ {\em is} singular at
Fuchsian groups on the boundary of $\cal P( N,\ua)$.

\medskip

The origin of these results was~\cite{KSQF}, which proved  the above
theorems in the very special case of  quasifuchsian once-punctured
tori, using much more elementary techniques. In~\cite{DS} we carried out
direct computations which proved Theorem~\ref{thm:globalpar} for some
very special curve systems on the twice punctured torus.
The results should have various applications. For example, combining
Theorem~\ref{thm:globalpar} with~\cite{S2}, when the holonomy of $N$
is quasifuchsian, one should be able to  exactly locate $\cal
P(N,\ua)$ in $\cal R(N)$.

As in~\cite{BonO}, our main tool is the local deformation theory of
cone manifolds developed by Hodgson and Kerckhoff in~\cite{HK}.  Let
$M$ be the $3$-manifold obtained by first doubling $\bar N$ across its
boundary and then removing the curves $\ua$.  It is easy to see that
a convex structure on $\bar N$ gives rise to a cone structure on $M$
with singular locus $\ua$. This means that everywhere in $M$ there are
local charts to $\HH^3$, except near a singular axis $\a$, where there
is a cone-like singularity with angle $2(\pi-\theta_{\a})$.  Under the
developing map, the holonomy of the meridian $m_{\a}$ around $\a$ is
an elliptic isometry with rotation angle $2(\pi-\theta_{\a})$. The
local parameterization theorem of Hodgson and Kerckhoff (see our
Theorem~\ref{thm:kh}) states that in a neighborhood of a cone
structure with singular axes $\ua=\{\a_1,\ldots,\a_d\}$ and cone
angles at most $2 \pi$, the representation space ${\cal R}(M)$ is
locally a smooth complex variety of dimension $d$, parameterized by
the complex lengths $\mu_{\a}$ of the meridians.  Notice that for a
cone structure, $\mu_{\a}$ is $\sqrt {-1}$ times the cone angle.
Moreover, the condition that the $\mu_{\a}$ are purely imaginary
characterizes the cone structures. This leads to a local version of
the Bonahon-Otal parameterization of convex structures in terms of
bending angles.

To prove Theorem~\ref{thm:main2}, we use the full force of the
holomorphic parameterization of $\R(M)$ in terms of the $\mu_{\a}$.
Under the hypothesis that $\ua=\{\a_1,\ldots,\a_d\}$ is maximal, the
spaces $\R(M)$ and $\R(N)$ have the same complex dimension $d$.
Moreover, we have the natural restriction map $r: \R(M) \rightarrow
\R(N)$.  Consider the pull-back $F=\cal L \circ r : \R(M) \rightarrow
\CC^d$ of the complex length function to $\R(M)$. Let $\rho_0 \in
\R(M)$ be the holonomy representation of the cone structure obtained
by doubling the convex structure $\sigma_0 \in \P(N,\ua)$.
Theorem~\ref{thm:main2} will follow if we show that $F$ is a local
diffeomorphism near $\rho_0$.  The key idea is that $F$ is a `real
map', that is, having identified the cone structures in $\R(M)$ with
$\RR^d$ (strictly speaking with $(i\RR)^d$), it has the properties:
\begin{eqnarray}
F(\RR^d) &\subset& \RR^d \label{eqn:easy}\\
F^{-1}(\RR^d) &\subset& \RR^d. \label{eqn:hard}
\end{eqnarray}
Using the fact that $F$ is holomorphic, we show in
Proposition~\ref{prop:holo} that this is sufficient to guarantee
that $F$ has no branch points and is thus a local diffeomorphism.

The first inclusion~(\ref{eqn:easy}) is relatively easy; the local
parameterization by cone angles allows us to show that, near the
double of a convex structure, the holonomy of any curve fixed by the
doubling map has real trace (Corollary~\ref{cor:realtrace}). Now
consider the inclusion~(\ref{eqn:hard}). We factor $F^{-1}(\RR^d)$ as
$r^{-1} ( \cal L^{-1}(\RR^d))$ and consider each of the preimages
separately.  In Section~\ref{sec:localpleating} we use geometrical
methods to prove the {\em local pleating theorem},
Theorem~\ref{thm:localpleating}, which states that near $\sigma_0$,
convex structures are characterized by the condition that the complex
lengths of the curves in the bending locus are real. (Actually we have
to introduce a slightly more general notion of a {\em piecewise
  geodesic structure}, in which we allow that some of the bending
angles may be negative, corresponding to some cone angles greater than
$2\pi$.) Thus $\L^{-1}(\RR^d)$ is locally contained in $\P(N,\ua)$.
The second main step in our proof, Theorem~\ref{thm:isomorphism}, uses
the duality between the meridians and the curves in the bending locus
to show that $r$ is a local holomorphic bijection between $\R(M)$ and
$\R(N)$. To finish the proof of (2), note that every convex structure
is the restriction of {\em some} cone structure, namely, its double.
However, since by Theorem~\ref{thm:isomorphism} the map $r$ is
one-to-one, the inverse image $r^{-1}(\sigma)$ of a convex structure
$\sigma$ near $\sigma_0$ can only be its double, and thus a cone
structure as required.

\medskip The plan of the rest of the paper is as follows.
Sections~\ref{sec:convexstr} and~\ref{sec:defspaces} supply
definitions and background. In Section~\ref{sec:convexity} we prove
the topological characterization Theorem~\ref{thm:nonemptypleating} of
which curve systems can occur, giving some extra details in the Appendix.
In order to simplify subsequent arguments, we also show
(Proposition~\ref{prop:lifting}) that the holonomy representation of a
cone structure coming from doubling a convex structure can always be
lifted to a representation into $SL(2,\CC)$.
Section~\ref{sec:defspaces} contains a brief review of the relevant
deformation theory, in particular expanding on the precise details of
the Hodgson-Kerckhoff parameterization near cusps.

The deduction of the
global parameterization theorem, Theorem~\ref{thm:globalpar}, from the
local version Theorem~\ref{thm:main2} is carried out in
Section~\ref{sec:main}.
The main idea is to observe that the Jacobian of the map $\cal L$
restricted to $\cal P(N,\ua)$, is the Hessian of the volume.
We deduce that it is positive definite and symmetric,
from which follow both the injectivity in
Theorems~\ref{thm:globalpar} and~\ref{cor:mixedparams}
and some additional information on volumes of convex cores.

   \medskip
\noindent{\bf Acknowledgments}  We would like to thank
Steve Kerckhoff and Max Forester for their insightful comments and
help, always dispensed with much generosity.  We would also like to
thank Darryl McCullough for detailed suggestions about the Appendix.
The second author is grateful for the support of her EPSRC Senior
Research Fellowship.


\section{Convex structures and their doubles}
\label{sec:convexstr}
\subsection{The bending locus}
\label{sec:convexcore}
The background material in this section is explained in detail in
~\cite{MT} and ~\cite{Mor}. We give a brief summary of the notions we
will need.  Let $G$ be a geometrically finite Kleinian group of the
second kind, so that its regular set $\Omega(G)$ is non-empty.  Let
$\cal H(G)$ be the hyperbolic convex hull of the limit set
$\Lambda(G)$ in $\HH^3$, and assume that $\Lambda(G)$ is not contained
in a circle, or equivalently, that the interior of $\cal H(G)$ is
non-empty.  The {\em convex core} of $N(G)=\HH^3/G$ is the
$3$-manifold with boundary $V=V(G) = \cal H(G)/G$. Alternatively,
$V(G)$ is the smallest closed convex subset of $N(G)$ which contains
all closed geodesics. As stated in the introduction, we assume that
$V(G)$ contains no rank-$2$ cusps.

As a consequence of the Ahlfors finiteness theorem, the boundary $\dd
V$ of $V$ consists of a finite union of surfaces, each of negative
Euler characteristic and each with possibly a finite number of
punctures.  Geometrically, each component $S_j$ of $\dd V$ is a
pleated surface whose bending locus is a geodesic lamination $pl(S_j)$
on $S_j$, see for example~\cite{EpM}.

In this paper we confine our attention to the case of {\em rational}
bending loci, in which the bending lamination of each $S_j$ is a set
of disjoint, homotopically distinct simple closed geodesics
$\{\a_{j_1}, \ldots, \a_{j_{n_j}}\}$.  We will use only the
union of these curves, renumbering them as $ \{\a_1, \ldots, \a_{m}
\}$, so that $m = \sum n_j$.  Let $\theta_i \in (0,\pi)$ denote the
exterior bending angle on $\a_i$ measured so that $\theta_i
\rightarrow 0$ as the outwardly oriented facets of $\dd V$ meeting
along $\a_i$ become coplanar.

It is convenient to modify the set $ \{\a_1, \ldots, \a_{m} \}$ so as
to include the rank-$1$ cusps of $\dd V$ as follows.  Let $\bar N(G)$
be the compact $3$-manifold with boundary obtained from $V$ by
removing a horoball neighborhood of each cusp.  Since $G$ is
geometrically finite, $\bar N(G)$ is compact; the notation has been
chosen to indicate that its interior is homeomorphic to $N(G)$.  More
precisely, choose $\epsilon>0$ sufficiently smaller than the Margulis
constant, so that the $\epsilon$-thin part $V_{(0,\epsilon]}$
(consisting of points in $V$ at which the injectivity radius is at
most $\epsilon$) consists only of finitely many disjoint rank-$1$
cusps.  Define $\bar N(G)$ to be the underlying manifold of the
$\epsilon$-thick part $V_{[\epsilon,\infty)}$. Note that for any
$\epsilon' < \epsilon$, there is a strong deformation retract of
$V_{[\epsilon',\infty)}$ onto $V_{[\epsilon,\infty)}$, and hence as a
topological manifold, $\bar N(G)$ is well-defined, independent of the
choice of $\epsilon$.

\begin{figure}[htb]
\centerline{\epsfbox{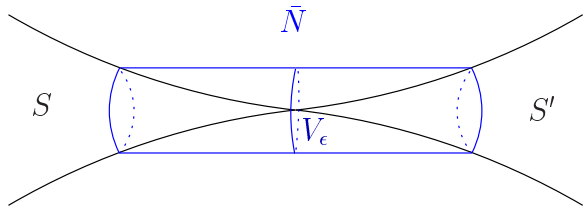}}
\caption{Replacing a pair of punctures with an annulus.}
\label{fig:annuli}
\end{figure}

The intersection $V_{\epsilon} = V_{[\epsilon,\infty)} \cap
V_{(0,\epsilon]}$ consists of the incompressible annuli which come
from the rank-$1$ cusps, see~\cite{Mor} for more details.  The
boundary $\dd \bar N(G)$ consists of the closed surfaces obtained from
$\dd V$ by removing horospherical neighborhoods of pairs of punctures
and replacing them with the annuli in $V_{\epsilon}$, as shown in
Figure~\ref{fig:annuli}. Conversely, we see that $\dd V$ can be
obtained from $\dd\bar N(G)$ by pinching the core curves of
$V_\epsilon$.
Denote the core curves in $V_\epsilon$ by
$\ua_P(G)=\{\a_{m+1},\ldots,\a_n\}$.

We define the {\em bending locus of $G$} to be the set 
$pl(G) =  \{\a_1,\ldots,\a_n\}$. 
We assign the bending angle  $\pi$ to each curve in $\ua_P(G)$,
corresponding to the fact that two facets of $\dd V$ which meet in a
rank-$1$ cusp as in Figure~\ref{fig:annuli}, lift to planes in
$\HH^3$ which are tangent at infinity, in which case the angle between
their outward normals is $\pi$.

\subsection{Definition of a convex structure}
\label{sec:defconvex}
In order to study variations of the structure defined in the last
section, we now make a more general topological discussion.  Let $\bar
N$ be a compact orientable $3$-manifold whose interior $N$ admits a
geometrically finite hyperbolic structure.  Assume that $\dd\bar N$ is
non-empty and contains no tori.  A {\em curve system} on $\dd\bar N$
is a collection of disjoint simple closed curves $\u{\a} = \{
\a_1,\ldots, \a_n\}$ on $\dd\bar N$, no two of which are homotopic in
$\dd\bar N$. In light of the previous section, we designate a subset
$\ua_P$ of $\ua$ as the {\em parabolic locus}.  The previous section
describes a generic {\em convex structure} on $(\bar N, \ua)$. We wish
however, to be slightly broader, by allowing the bending angles on
some of the curves in $\ua$ to vanish. It is convenient to begin with
the following slightly more general definition. Note that both $\bar
N$ and $\bar N - \ua_P$ are homotopically equivalent to $N$.
\begin{defn} 
\label{def:piecewisegeo}
A  piecewise geodesic structure on $(\bar N,\u{\a})$ 
with parabolic locus $\ua_P$ consists of a pair $(G, \phi)$, where $G$
is a geometrically finite Kleinian group whose convex core has
non-zero volume and where
$\phi: \bar N - \ua_P  \hookrightarrow N(G)$ is an embedding
such that the following properties are satisfied:\\
(i) $\phi_*:\pi_1(N) \rightarrow \pi_1(N(G))$ is an
isomorphism; \\
(ii)  $\phi_*(\gamma)$ generates a maximal parabolic subgroup for
every $\gamma$ freely homotopic to a \hspace*{0.6cm} curve in $\ua_P$;\\ 
(iii) the image $\phi(\a)$ of each $\a \in \ua - \ua_P$ is geodesic
and the image  of each component 
\hspace*{0.6cm} of $\dd\bar N- \ua$ is totally geodesic.
\end{defn}

The closures in $\HH^3$ of the lifts of the components of $\phi(\dd
\bar N- \ua)$ to $\HH^3$ will be called {\em plaques}. Clearly, each
plaque is totally geodesic and each component of its boundary projects
to a geodesic $\phi(\a), \a \in \ua - \ua_P$. Notice that the plaques
are not necessarily contained in $\dd {\cal H}(G)$. However, since
$\phi $ is an embedding, two plaques can only intersect along a common
boundary geodesic. We call such  geodesics and their projections 
{\em bending lines}.

Let $x$ be a point on a bending line and let $B$ be a small
ball containing $x$. The two plaques form a `roof' which separates
$B$ into two components, exactly one of which is contained in Image$(\phi)$. 
Let $\psi_\a \in (0,2\pi)$ be the angle between the plaques on
the side intersecting Image$(\phi)$ and let $\theta_{\a}=\pi -
\psi_\a$, so that $\theta_{\a} = 0$ when the plaques are coplanar. 
We call $\theta_{\a}$ the {\em bending angle} along $\a$. Notice it
is possible that $\theta_{\a} < 0$. If $\a \in \ua_P$ we again set 
$\theta_{\a} = \pi$.

We will often allude to a piecewise geodesic structure on $(\bar
N,\ua)$ without mentioning its parabolic locus.  This simply means
that it is a piecewise geodesic structure on $(\bar N,\ua)$ with some
parabolic locus, which is unnecessary to specify.  A piecewise
geodesic structure determines a holonomy representation $\pi_1(N)
\rightarrow PSL(2,\CC)$ up to conjugation.  Two piecewise geodesic
structures are equivalent if their holonomy representations are
conjugate.

A piecewise geodesic structure is a generalization of 
a convex structure. In terms of the above definition, we have:
\begin{defn}
  A convex structure $(G,\phi)$ on $(\bar N,\u{\a})$ with parabolic
  locus $\ua_P$ is a piecewise geodesic structure on $(\bar N,\ua)$
  with the same parabolic locus, which satisfies the additional
  property that the image of $\phi$ is convex.
\end{defn} 

Convexity is easily described in terms of bending angles; since local
convexity implies convexity, see for example~\cite{CEG} Corollary
1.3.7, a piecewise geodesic structure is convex if and only if
$\theta_{\a} \geq 0$ for all $\a \in \ua$.  It follows easily from
the definition that if $(G,\phi)$ is a convex structure on $(\bar
N,\ua)$, then the image of $\phi$ is equal to the convex core $V(G)$.
Furthermore, $pl(G)$ will be contained in $\ua$ and the core curves
$\ua_P(G)$ of the annuli around the rank-$1$ cusps described in the
previous section, is the parabolic locus of the convex structure.  In
other words, we are exactly in the situation described in
Section~\ref{sec:convexcore}, except that we have allowed ourselves to
adjoin to the bending locus some extra curves along which the bending
angle is zero.

\subsection{Topological characterization of $(\bar N,\ua)$}
\label{subsec:topconds} 
As in the previous section,
let $\bar N$ be a compact orientable $3$-manifold whose interior $N$
admits a complete hyperbolic structure. We assume that $\dd \bar N$ is
non-empty and contains no tori.  Theorem~\ref{thm:nonemptypleating}
below gives a topological characterization of the curve systems
$\u{\a}$ on $\dd \bar N$ for which there exists some convex structure
on $(\bar N,\u{\a})$. Our statement is essentially a reformulation of
results in~\cite{BonO}.

First recall some topological definitions. 
A surface $F \subset \bar N$ is {\em properly embedded} if $\dd F
\subset \dd \bar N$ and $F$ is transverse to $ \dd \bar N$.
An {\em essential disk} $D \subset \bar N$ is a properly embedded  disk 
which cannot be homotoped to a disk in 
$\dd \bar N$ by a homotopy fixing  $\dd D$. 
An {\em essential  annulus}  $A \subset \bar N$ is a   properly
embedded annulus which is not null homotopic in
$\bar N$, and  which cannot be homotoped into an annulus in 
$\dd \bar N$ by a homotopy which fixes $\dd A$.

\begin{defn}(Thurston,~\cite{Thu3})
Let $\bar N$ be defined as above.
A curve system  $\u{\a}$ on 
$\dd\bar N$ is {\em doubly incompressible}  with respect to $(\bar
N,\dd\bar N)$  if:\\
D.1 There are no essential annuli with boundary in
$\dd\bar N -\u{\a}$.\\
D.2 The boundary of every essential disk intersects $\u{\a}$ at least
 $3$ times.
\end{defn}
The characterization is as follows:
\begin{thm}
\label{thm:nonemptypleating} Let $\bar N$ be defined as above and let
$\u{\a}$ be a non-empty curve system on
$\dd\bar N$. There is a convex structure on 
$(\bar N, \u{\a})$ if and only if $\u{\a}$  is
doubly incompressible with respect to $(\bar N,\dd\bar N)$.
\end{thm}
{\bf Remark} Thurston's original definition has a third condition
(D.3) which states that every maximal abelian subgroup of $\pi_1(\dd
\bar N -\ua)$ is mapped to a maximal abelian subgroup of $\pi_1(\bar
N)$. However, it can be shown that (D.1) implies (D.3). We thank the
referee for pointing this out.

\medskip

The remainder of this subsection outlines a proof of 
Theorem~\ref{thm:nonemptypleating}.

The necessity of the condition on $\u{\a}$ is a consequence of the
following result proved in~\cite{BonO}, whose proof we briefly
summarize for the reader's convenience:
\begin{prop}
\label{thm:BOconditions} (\cite{BonO} Propositions 4 and 7.) 
Let $\ua = pl(G)$ be the bending lamination of a geometrically finite
Kleinian group $G$, and let $\theta_i \in (0,\pi]$ be the bending
angle on $\a_i$. Let $\bar N(G)$ be the associated compact
$3$-manifold as defined in Section~\ref{sec:convexcore}, and let $\xi
$ be the measured lamination which assigns the weight $\theta_i$ to
each intersection with the free
homotopy class of the curve $\a_i$ on $\dd \bar N(G)$. Then:\\
(i) For each essential annulus $A$ in $\bar N(G)$, we have $i(\dd A,
\xi)>0$.\\
(ii) For each essential disk $D$ in $\bar N(G)$, we have $i(\dd D,
\xi) >2 \pi$.
 \end{prop}
 Here $i(\gamma,\xi)$ denotes the intersection number of a loop
 $\gamma$ with the measured lamination $\xi = \sum_i \theta_i \a_i$.
 We remark that in \cite{BonO}, the theorem does not assume that the
 convex core $V(G)$ contains no rank-$2$ cusps. In that case, $\bar
 N(G)$ is defined analogously as the manifold obtained from $V(G)$ by
 removing disjoint horoball neighborhoods of both rank-$1$ and
 rank-$2$ cusps.  \medskip

\noindent {\sc Sketch of Proof.} If $i(\dd A,
\xi)=0$ then the two components of $\dd A$ are either freely homotopic
to geodesics in $\HH^3/G$ or loops around punctures of $\dd V(G)$. It
is impossible for one component of $\dd A$ to be geodesic and one to
be parabolic. If both components are geodesic, lifting to $\HH^3$, we
obtain an infinite annulus whose boundary curves are homotopic
geodesics at a bounded distance apart, which therefore coincide. The
resulting infinite cylinder bounds a solid torus in $\HH^3$, from
which one obtains a homotopy of $A$ into $\dd \bar N(G)$, showing that
$A$ was not essential.  Finally, if both components of $\dd A$ are
parabolic, they can only be paired in a single rank-$1$ cusp of $V(G)$
from which it follows that $A$ was not essential.

In the case of a disk, note that since $G$ is torsion free, $\dd D$ is
necessarily indivisible and moreover, cannot be a loop round a
puncture.  Thus if $i(\dd D, \xi)=0$, then $\dd D$ would be freely
homotopic to a geodesic in $\HH^3/G$, which is impossible. Now
homotope $\dd D$ to be geodesic with respect to the induced hyperbolic
metric on $\dd \bar N$.  We can also homotope $D$ fixing the boundary
so that it is a pleated disk which is a union of totally geodesic
triangles. Note that the interior angles between two consecutive
segments in $\partial D$ are greater than the dihedral angles between
the corresponding planes.  The Gauss Bonnet Theorem applied to the
pleated disk now gives the result. \ q.e.d. \vskip .2in

\begin{cor}
\label{cor:necnonemptypleating} Let $\bar N$ be defined as above
and let $\u{\a}$ be a curve system on
$\dd\bar N$. If there exists a convex structure  on $(\bar N,
\u{\a})$, then  $\u{\a}$  is
doubly incompressible with respect to $(\bar N,\dd\bar N)$.
\end{cor}
\begin{proof}  We check the
conditions for $\u{\a}$ to be doubly incompressible.
The bending measure of each curve is at most $\pi$. Thus ({\em D.1}\,)
and ({\em D.2}\,) follow from conditions (i) and (ii) of
Proposition~\ref{thm:BOconditions}, respectively. Strictly speaking,
conditions (i) and (ii) imply that ({\em D.1}\,) and ({\em D.2}\,) hold for 
the curve system $\ua' \subset \ua$ on which $\theta_{\a'}>0$.
However, if ({\em D.1}\,),({\em D.2}\,) hold for the subset $\ua'$, then they
certainly hold for the larger set $\ua$. 
\end{proof}
Conversely, we have:
\begin{prop}
\label{prop:sufnonemptypleating} Let $\bar N$ be as defined above.
If $\u{\a}$ is a doubly incompressible curve system with respect to
$(\bar N,\dd\bar N)$, then there is a convex structure on $(\bar N,
\u{\a})$ for which $\ua=\ua_P$.
\end{prop}
The idea of the proof is to show that the conditions on $\ua$
guarantee that the manifold $M$ obtained by first doubling $\bar N$
across $\dd \bar N$ and then removing $\ua$ is both irreducible and
atoroidal. It then follows from Thurston's hyperbolization theorem for
Haken manifolds that $M$ admits a complete hyperbolic structure. It is
not hard to show that this structure on $M$ induces the desired convex
structure on $(\bar N,\ua)$. The essentials of the proof are contained
in~\cite{BonO} Th{\'e}or{\`e}me 24.  For convenience we repeat it in
the Appendix, at the same time filling in more topological details.

\label{sec:convexity}
\subsection{Doubles of Convex Structures}
\label{sec:doubles} 
A convex structure on $(\bar N,\ua)$ naturally induces a {\em cone
  structure} on its double. Topologically, we form the {\em double} $D
\bar N$ of $\bar N$ by gluing $\bar N$ to its mirror image $\tau(\bar
N)$ along $\dd\bar N$.  We may regard $\tau$ as an orientation
reversing involution of $D\bar N$ which maps $\bar N$ to $\tau(\bar
N)$ and fixes $\dd\bar N$ pointwise.

A convex structure on $(\bar N, \u{\a})$ clearly induces an isometric
structure 
on $\tau (\bar N)$. Since gluing $\bar N$ and $\tau(\bar N)$
along $\dd\bar N$ matches the hyperbolic structures everywhere except
at points in $ \u{\a}$, this naturally induces a {\em cone structure}
on $M=D \bar N -\ua$.  If the bending angle along $\a \in \ua$ is
$\theta_\a$, then the cone angle $\varphi_\a$ around $\a$ is
$2(\pi-\theta_\a)$.
More precisely, a {\em hyperbolic cone
  structure} on $M$ with singular locus $\ua$ is an incomplete
hyperbolic structure on $M$ whose metric completion determines a
singular metric on $D\bar N$ with singularities along $\ua$.  Often we
refer to this simply as a cone structure on $(M,\ua)$.  In the
completion, each loop $\a \in \ua$ is geodesic and in cylindrical
coordinates around $\a$, the metric has the form
\begin{equation}
\label{eqn:conemetric} 
dr^2 + \sh ^2 r d\t^2 + \ch ^2 r dz^2, \end{equation}
where $z$ is the distance along the
singular locus $\a$, $r$ is the distance from $\a$, 
 and $\theta$ is the angle around $\a$
measured modulo some $\varphi_\a > 0$. The angle $\varphi_\a$ is  called
the {\em cone angle} along $\a$, see~\cite{HK}. More generally,
a cone structure is allowed have cone angle zero along 
a subset of curves in $\ua$, which in our case will be the parabolic
locus $\ua_P$ of the convex structure. 
This means that the metric completion determines a singular
metric on the interior of $D \bar N -\ua_P$ with singularities along
$\ua-\ua_P$ as described in Equation~(\ref{eqn:conemetric}). 
The metric in a neighborhood of a missing curve
$ \a  \in \ua_P$ is complete, making it a rank-$2$ cusp.

Geometric doubling can be just as easily carried out for a piecewise
geodesic structure on $(\bar N,\ua)$. The only difference is that the
resulting cone manifold may have some cone angles greater than $2\pi$.

Associated to a cone structure is a developing map $dev: \widetilde{M}
\to \HH^3$ and a holonomy representation $\rho: \pi_1(M) \to
PSL(2,\CC)$.  It is well known that if the image of a representation
$\sigma :\pi_1(\bar N) \to PSL(2,\CC)$ is torsion free and discrete,
then it can be lifted to $SL(2,\CC)$, see for example~\cite{Cu}.  In
the remainder of this section, we prove that, even though in general
it is neither free nor discrete, the holonomy representation of a cone
manifold $M$ formed by geometric doubling can also be lifted. This
conveniently resolves any difficulties about defining the trace of
elements $\rho(\gamma)$ later on.  In the course of the proof, we
shall find an explicit presentation for $\pi_1(M)$ in terms of
$\pi_1(\bar N)$, and then describe explicitly how to construct the
holonomy $\rho$ of $M$ starting from the holonomy $\sigma$ of $\bar
N$.

In general, suppose that $G=\langle g_1,\ldots,g_k\ |\ R_1,\ldots,R_l
\rangle$ is a finitely presented group.  To lift a homomorphism
$\phi:G\to H$ to a covering group $\tilde H \to H$, we have to show
that for each generator $g_i$ we can choose a lift $\tilde \phi(g_i)
\in \tilde H$ of $\phi(g_i)$ in such a way that $\tilde \phi(g_{i_1} )
\ldots \tilde \phi(g_{i_s} ) = {\rm id}_{\tilde H}$ for each relation
$R_i= g_{i_1} \ldots g_{i_s} = {\rm id} $ in $G$.

\begin{prop}
\label{prop:lifting}
Suppose that the holonomy representation $\sigma$ of a convex
structure on $(\bar N,\u{\a})$ lifts to a representation $\tilde
\sigma: \pi_1(\bar N) \to SL(2,\CC)$. Then the holonomy representation
$\rho$ for the induced cone structure on its double also lifts to a
representation $\tilde \rho: \pi_1( M) \to SL(2,\CC)$.
\end{prop}
\begin{proof} We begin by finding an explicit presentation for 
  $\pi_1(M)$ in terms of $\pi_1(\bar N)$.  Let the components of $\dd
  \bar N - \u{\a}$ be $S_0,\ldots, S_k$.  We will build up $\pi_1(M)$
  first by an amalgamated product and then by HNN-extensions by
  glueing $\bar N-\ua$ to its mirror image $\tau(\bar N-\ua)$ in
  stages, at the $i^{ th}$ stage gluing $S_i$ to $\tau(S_i)$.
  
  For each $i$, pick a base point $x_i \in S_i$ and pick paths
  $\beta_i$ from $x_0$ to $x_i$ in $\bar N-\ua$.  Then $\tau(\beta_i)$
  is a path from $\tau(x_0)$ to $\tau(x_i)$ in $\tau(\bar N)$.  Let
  $\tau_i = \tau|_{S_i}$.  First, glue $S_0$ to $\tau(S_0)$ using
  $\tau_0$ to form a manifold $M_0$. Then $\pi_1(M_0, x_0)= \pi_1(\bar
  N, x_0) *_{\pi_1(S_0,x_0)} \pi_1(\tau(\bar N), \tau(x_0))$, where
  $\pi_1(S_0) \to \pi_1(\bar N)$ is induced by the inclusion $\iota:
  S_0 \hookrightarrow \bar N$ and $\pi_1(S_0) \to \pi_1(\tau(\bar N))$
  is induced by $\tau \circ \iota$.
  
  Now suppose inductively we have glued $S_{i-1}$ to $\tau(S_{i-1})$
  forming a manifold $M_{i-1}$ and that we know $\pi_1(M_{i-1},x_0)$.
  Now glue $S_i$ to $\tau(S_i)$ using $\tau_i$ and denote the
  resulting manifold $M_i$.  This introduces a new generator $e_i =
  \b_i \tau(\b_i)^{-1}$.  Let $G_i$ denote the image of
  $\pi_1(S_i,x_i)$ in $\pi_1(M_{i-1},x_0)$ under the inclusion map,
  where loops based at $x_i$ are mapped to loops based at $x_0$ by
  concatenating with $\beta_i$. Then van-Kampen's theorem implies that
  the fundamental group $\pi_1(M_i,x_0)$ has the presentation $\langle
  \pi_1(M_{i-1},x_0) , e_i\,|\, e_i^{-1} \gamma e_i=\tau(\gamma),
  \gamma \in G_i \rangle$.  In this way, we inductively obtain a
  presentation for $\pi_1(M,x_0)$.

\smallskip

We now want to give an explicit description of the holonomy
representation $\rho$ for the doubled cone structure on $M$ in terms
of the holonomy representation $\sigma$ for the convex structure on
$\bar N$.  First consider $\sigma$. The base point $x_0$ of $\bar N$
is contained in a totally geodesic plaque $\Pi$ in the convex core
boundary.  The developing map $dev:\widetilde{ \bar N} \to \HH^3$ and
resulting holonomy representation $\sigma:\pi_1(\bar N,x_0)\rightarrow
PSL(2,\CC)$ are completely determined by a choice of the image of
$dev(x_0)$ and image of an inward pointing unit normal $\underline n$
to $\Pi$ at $x_0$. Let $dev_{\tau}$ be the developing map of
$\tau(\bar N)$ for which $dev_{\tau} \tau(x_0) = dev(x_0)$ and
$dev_{\tau}(\underline n) = -dev(\underline n)$.  Now $ dev(x_0)$ lies
in a hyperbolic plane which is fixed by the Fuchsian subgroup
$\sigma(\pi_1(G_0))$. Let $J$ be inversion in this plane. Then
$J(dev(\underline n)) = -dev(\underline n)$ and we deduce that
$dev_{\tau} \circ \tau = J \circ dev$ and hence the associated
holonomy representation $\hat \sigma: \pi_1(\tau(\bar N)) \rightarrow
PSL(2,\CC)$ is given by $\hat \sigma \circ \tau_*=J\sigma J^{-1}$.

Clearly, $dev$ and $dev_{\tau}$ together determine $\rho$. Our
explicit description of $\rho$ will be found by inductively finding
the holonomy representation $\rho_i$ for the induced cone structure on
$M_i$, $i=0,1,\ldots,k$.

Define a representation $ \pi_1(\bar N, x_0) *\pi_1(\tau(\bar N
),\tau(x_0)) \to PSL(2,\CC)$ by specifying that its restrictions to $
\pi_1(\bar N,x_0)$ and $\pi_1(\tau(\bar N),\tau(x_0))$ are $\sigma$
and $\hat \sigma$ respectively.  Since $J(\gamma) = \gamma$ for
$\gamma \in \pi_1(S_0,x_0)$, we deduce from the amalgamated product
description of $\pi_1(M_0,x_0)$ above that $\sigma * \hat \sigma$
descends to a representation $\rho_0:\pi_1(M_0,x_0) \rightarrow
PSL(2,\CC)$. This is clearly the holonomy representation of $M_0$.
Now, suppose inductively that we have found the holonomy
representation $\rho_{i-1}:\pi_1(M_{i-1},x_0) \rightarrow PSL(2,\CC)$.
From the HNN-extension description of $\pi_1(M_{i},x_0)$, we see that
in order to compute $\rho_{i}$, it is sufficient to find $\rho(e_i)$.
It is not hard to see that $\rho(e_i)= J_i J$, where $J_i$ is the
orientation reflection in the plane through $dev(x_i)$ which is fixed
by the Fuchsian group $\sigma (G_i)$. The holonomy representation
$\rho$ of $\pi_1(M)$ is equal to $\rho_k$ found inductively in this
way.

\smallskip

Finally, this careful description allows an easy solution of the
lifting problem. Given a lifting $\tilde \sigma: \pi_1( \bar N) \to
SL(2,\CC)$ of the holonomy representation $\sigma: \pi_1( \bar N) \to
PSL(2,\CC)$, we want to define a corresponding lifted representation
$\tilde \rho$ of $\rho$.  Following the inductive procedure for
constructing $\rho$ above, we see that the only requirement on
$\rho(e_i)$ is that it satisfy the relation
$\rho(e_i)^{-1}\rho(\gamma)\rho(e_i) = \rho(\tau(\gamma))$ for all
$\gamma \in G_i$. Thus we have to show that at each stage the isometry
$\rho(e_i) = J_i J$ can be lifted to an element $\tilde \rho(e_i) \in
SL(2,\CC)$ which satisfies the relation
$$\tilde \rho(e_i)^{-1}\tilde\rho(\gamma)\tilde\rho(e_i)=
\tilde\rho(\tau(\gamma)).$$
Since $\rho(\tau(\gamma)) = J \rho(\gamma)
J$, this relation reduces in $PSL_2(\CC)$ to $J_i$ commuting with
$\rho(\gamma) = \sigma(\gamma)$ for all $\gamma \in G_i$.  This just
means that $J_i$ fixes axes of elements in $\sigma (G_i)$, which is
clearly the case.  The lifted relation is obviously satisfied
independently of the choice of $\tilde \sigma(\gamma)$ and for either
choice of lift of $J_iJ$, which is all we need.
\end{proof}

The following general fact about lifting is also clear:
\begin{prop} Suppose that a representation $\rho_0 \in Hom(\pi_1(M),
  PSL(2,\CC) )$ lifts to $ SL(2,\CC)$. Then $\rho_0$ has a
  neighborhood in $\mbox{Hom}(\pi_1(M), PSL(2,\CC))$ in which every
  representation also lifts to $SL(2,\CC)$.
\end{prop}
Here, $Hom(\pi_1(M),PSL(2,\CC))$ is the space of homomorphisms from
$\pi_1(M)$ to $SL(2,\CC)$. It has the structure of a complex variety,
which is naturally induced from the complex structure on $SL(2,\CC)$.


\section{Deformation spaces}
\label{sec:defspaces}
Let $\bar N$ be a compact orientable $3$-manifold whose interior $N$
admits a complete hyperbolic structure. As usual, assume that $\dd\bar
N$ is non-empty and contains no tori.  Let $\ua$ be a doubly
incompressible curve system on $\dd\bar N$ and let $M=D\bar N -\ua$.
The possible hyperbolic structures on $N$ and cone structures on $M$
are locally parameterized by their holonomy representations
$\sigma:\pi_1(N) \rightarrow PSL(2,\CC)$ and $\rho:\pi_1(M)
\rightarrow PSL(2,\CC)$ modulo conjugation. As shown in
Proposition~\ref{prop:lifting}, all the representations relevant to
our discussion can be lifted to $SL(2,\CC)$.  Thus from now on, to
simplify notation, we shall use $\sigma$ and $\rho$ to denote the {\em
  lifts} of the holonomy representations of $N, M$ respectively to
$SL(2,\CC)$.

Let $W$ denote either $N$ or $M$ and consider the space of
representations $\R(W) = \mbox{Hom}(\pi_1(W),SL(2,\CC))/ SL(2,\CC),$
where $SL(2,\CC)$ acts by conjugation.  When it is necessary to make a
distinction, the equivalence class of $\omega \in \mbox{Hom}(\pi_1(W),
SL(2,\CC))$ will be denoted $[\omega]$, although to simplify notation
we often simply write $\omega \in \R(W)$. Although in general, $\R(W)$
may not even be Hausdorff, in the cases of interest to us the results
below show that it is a smooth complex manifold.  (Section 3 of the
survey \cite{goldmanrep} is a good reference for further details.)

First we consider $\cal R(N)$.  If $\ua_P$ is a fixed curve system on
$\dd \bar N$, let $P$ be the set of elements in $\pi_1(N)$ which are
freely homotopic to a curve in $\ua_P$. We denote by $\R_P(N)$ the
image in $\cal R(N)$ of the set of all representations $\sigma:
\pi_1(N) \to SL(2,\CC)$ for which $\sigma(\gamma)$ is parabolic for
all $\gamma$ in $P$.

Now let $\sigma: \pi_1(N) \rightarrow SL(2,\CC)$ be the holonomy
representation for a geometrically finite structure on $N$. Put
$G=\sigma(\pi_1(N))$ and let $\ua_P$ be the collection of core curves
of the annuli in the parabolic locus of $V(G)$. From the Marden
Isomorphism Theorem~\cite{Marden}, we have that a neighborhood of
$[\sigma]$ in $\R_P(N)$ can be locally identified with the space of
quasiconformal deformations of $N(G)=\HH^3/G$. By Bers' Simultaneous
Uniformization Theorem ~\cite{Bers2}, this space is isomorphic to
$\Pi_i \teich(S_i)/ Mod_0 (S_i)$, where $S_i$ are the components of
$\dd V$, $\teich(S_i)$ is the Teichm\"uller space of $S_i$, and $Mod_0
(S_i)$ is the set of isotopy classes of diffeomorphisms of $S_i$ which
induce the identity on the image of $\pi_1(S_i)$ in $ \pi_1(N)$.  If
$S_i$ has genus $g_i$ with $b_i$ punctures, then $\teich(S_i)$ has
complex dimension $d_i = 3g_i-3 +b_i$. Note that $d_i$ is also the
maximal number of elements in a curve system on $S_i$. Thus we obtain:

\begin{thm} Let $N$ be a geometrically finite hyperbolic
  $3$-manifold with holonomy representation $\sigma: \pi_1(N) \to
  SL(2,\CC)$. Let $P$ be the set of elements in $\gamma \in \pi_1(N)$
  such that $\sigma(\gamma)$ is parabolic. Then $\R_{P}(N)$ is a
  smooth complex manifold near $[\sigma]$, of complex dimension
  $\sum_i d_i.$
  \end{thm}
  
  Since we wish to allow deformations which `open cusps', we will also
  need the smoothness of $\R(N)$. Since $\sigma(\pi_1(N))$ is
  geometrically finite, the punctures on the surfaces $S_i$ are all of
  rank-$1$ and they are all matched in pairs, see~\cite{Marden} and
  the discussion accompanying Figure~\ref{fig:annuli}. Thus the number
  of rank-$1$ cusps is $\sum_i b_i/2$ and opening up each pair
  contributes one complex dimension. Note that $\sum_i (d_i +b_i/2 )$
  is the number of elements in a maximal curve system on $\dd \bar N$.
  The following result is~\cite{Kap} Theorem 8.44, see also~\cite{Hod}
  Chapter 3:
\begin{thm}
  \label{thm:smoothreps} Let $N$ be a geometrically finite hyperbolic
  $3$-manifold with holonomy representation $\sigma: \pi_1(N) \to
  SL(2,\CC)$.  Then $\R (N)$ is a smooth complex manifold near
  $[\sigma]$, of complex dimension $\sum_i (d_i + b_i/2)$.
  \end{thm}
The special case in which $\pi_1(N)$ is a surface group
(so $\sigma$ is quasifuchsian) is treated in more detail in~\cite{Gol}.
We remark that in~\cite{Kap}, the above theorem is also stated in the
case in which $\dd \bar N$ contains tori.

\medskip

Now let us turn to the deformation space $\R(M)$ where $M$ is a cone
manifold as above.  The analogous statement to
Theorem~\ref{thm:smoothreps} in a neighborhood of a cone structure is
one of the main results in~\cite{HK}. In fact, Hodgson and Kerckhoff
give a local parameterization of $\R(M)$ by the complex lengths of the
{\em meridians}.  In terms of the coordinates in Equation
(\ref{eqn:conemetric}) in Section~\ref{sec:doubles}, a meridian
$m=m_{\a}$ is a loop around a singular component $\a$, which can be
parameterized as $(r(t),\theta(t),z(t))= (r_0,t,z_0)$ where $t \in
[0,\varphi_\a]$.  By fixing an orientation on $\a$, the meridian can
be chosen so that $m$ is a right-hand screw with respect to $\a$. To
define an element in $\pi_1(M,x_0)$, simply choose a loop in
$\pi_1(M,x_0)$ freely homotopic to $m$.  The particular choice is not
important, since we shall mainly be concerned with the complex length
or trace.  By abuse of notation, we shall often write $\rho_0(\a)$,
$\rho_0(m)$ to denote $\rho_0(\gamma)$ where $\gamma \in \pi_1(M)$ is
freely homotopic to $\a$ or $m$, as the case may be. We assume that
all representations concerned can be lifted to $SL(2,\CC)$.
\begin{thm}
\label{thm:kh} (\cite{HK} Theorem 4.7) Let $M$ be a finite volume
$3$-dimensional hyperbolic cone manifold whose singular locus is a
collection of disjoint simple closed curves $\u{\a} =
\{\a_1,\ldots,\a_n\}$.  Let $\rho_0: \pi_1(M) \rightarrow SL(2,\CC)$
be a lift of the holonomy representation for $M$.  If all the cone
angles $\varphi_i$ satisfy $0\leq \varphi_i \leq 2\pi$, then $\R(M)$
is a smooth complex manifold near $[\rho_0]$ of complex dimension $n$.
Further, if $m_1,\ldots, m_n$ are homotopy classes of meridian curves
and if $0<\varphi_i \leq2 \pi$, then the complex length map ${\cal M}:
\R(M) \rightarrow \CC^n$ defined by ${\cal M}([\rho])
=(\l(\rho({m_1})),\ldots,\l(\rho({m_n})))$ is a local diffeomorphism
near $[\rho_0]$.
\end{thm}
Here $\l(\rho({m_i}))$ denotes the complex length of $\rho(m_i)$,
discussed in more detail in the next section.  Structures for which
$\varphi_i=0$ are excluded from the local parameterization given by
the map $\cal M$ because strictly speaking, the complex length of
$\rho(m_i)$ cannot be defined as a holomorphic function in a
neighborhood of $[\rho_0]$ when $\rho_0(m_i)$ is parabolic.  However,
in such a case, replacing the complex length $\l(\rho({m_i}))$ with
its trace Tr$\rho(m_i)$ again gives a local parameterization of
$\R(M)$.  This and the case in which $\varphi_i=2\pi$ are expanded
upon in the next section (see also the proof of Theorem 4.5
in~\cite{HK} and the remark at the end of their section 4).

\subsection{Local deformations and complex length}
\label{subsec:complexlength}
  
The complex length $\l(A)$ of $A \in SL(2,\CC)$ is determined from its
trace by the equation $\tr A= 2 \ch \l(A)/2$. Since $ z \mapsto \cosh
z$ is a local holomorphic bijection except at its critical values
where $\cosh z = \pm 1$, the function $[\rho] \mapsto
\l(\rho(\gamma))$ is locally well-defined and holomorphic on the
representation space $\R(M)$, except possibly at points for which
$\rho(\gamma)$ is either parabolic or the identity in $SL(2,\CC)$.
  
If $\tr A \neq \pm 2$, then $\mbox{Re}\, \l(A)$ is the translation
distance of $A$ along its axis and $\mbox{Im}\, \l(A)$ is the
rotation.  The sign of both these quantities depends on a choice of
orientation for $\Ax A$, corresponding to the ambiguity in choice of
sign for $\l(A)$ in its defining equation.  For a detailed discussion
of the geometrical definition, see~\cite{Fenchel}~V.3 or~\cite{SWol}.

To study local deformations, we work at a point $[\rho_0] \in \R(M)$,
and study the possible conjugacy classes of one parameter families of
holomorphic deformations $t \mapsto [\rho_t] \in \R(M)$, defined for
$t$ in a neighborhood of $0$ in $\CC$.  For each $\gamma \in
\pi_1(M)$, the derivative $\dot{\rho}(\gamma)=\frac{d}{dt}\big|_{t=0}
(\rho_t(\gamma)\rho_0^{-1}(\gamma))$ is an element of the Lie algebra
$\sl2c$. In this way, an infinitesimal deformation defines a function
$\dot \rho =z:\pi_1(M) \rightarrow \sl2c$. The fact that
$\rho_t(\gamma_1 \gamma_2)=\rho_t(\gamma_1)\rho_t(\gamma_2)$ for all
$t$ and $\gamma_1, \gamma_2 \in \pi_1(M)$, forces $z$ to satisfy the
cocycle condition $z(\gamma_1 \gamma_2)=z(\gamma_1) +
Ad\r(\gamma_1)z(\gamma_2)$. The fact that holomorphically conjugate
representations are equal in $\R(M)$ implies that an infinitesimal
deformation with a cocycle of the form $z(\gamma)=v - Ad\r(\gamma)v$
for some $v \in \sl2c$, is trivial. Thus the space of infinitesimal
holomorphic deformations of $[\rho_0]$ in $\R(M)$ is identified with
the cohomology group $H^1(M;Ad \rho_0)$ of cocycles modulo
coboundaries. Moreover, if $\R(M)$ is smooth at $[\rho_0]$, then
$H^1(M;Ad \rho_0)$ can be identified with the holomorphic tangent
space $T_{[\rho_0]}\R(M)$.  A good summary of this material can be
found in~\cite{Kap}, see also~\cite{Gol} and~\cite{Hod}.
  
\smallskip Let us look in more detail at the parameterization of
$\R(M)$ given in Theorem~\ref{thm:kh}.  In~\cite{HK}, Corollary 1.2
combined with Theorems 4.4 and 4.5 shows that if $\rho$ is the
holonomy representation of a cone-structure with cone angles at most
$2\pi$, then $\Hom(\pi_1(M),SL(2,\CC))$ is a smooth manifold of
dimension $n+3$ near $\rho$ and that the restriction map $res:
H^1(M;Ad \rho) \to \oplus_{i=1}^n H^1(m_i;Ad \rho)$ is injective.
Specifying a parameterization is then only a matter of choosing a map
$\Phi:\R(M) \rightarrow \CC^n$ whose derivative can be identified with
$res$.  Expanding on the discussion in~\cite{HK}, we will verify that
the map $\Phi$ can be taken to be $\cal M$ defined above.  We
emphasize that if a cone angle $\varphi_i$ vanishes, then the
corresponding parameter should be changed from complex length
$\l(\rho(m_i))$ to trace $\tr \rho(m_i)$.

\smallskip

Denote the basis vectors of $ \sl2c$ as follows:
$$u^+ = \pmatrix{0 & 1 \cr 0 & 0}, \ \ u^-= \pmatrix{0 & 0 \cr 1 & 0},
\ \ v=\pmatrix{1 & 0 \cr 0 & -1}.$$
Also let $T_{\a_i}$ denote the
boundary torus of a tubular neighborhood of the singular axis $\a_i$.
For simplicity, in what follows, we drop the subscript $i$ so that
$\pi_1(T_\alpha)$ is generated by a longitude $\a$ and a meridian $m$.
As above, we use $\rho_0(\a)$, $\rho_0(m)$ to denote the image
$\rho_0(g)$ where $g \in \pi_1(M)$ is in the appropriate free homotopy
class.  Note that since $\rho_0$ will always be the double of a convex
structure, we may assume that $\rho_0(\a) \neq {\rm id}$.  \medskip

\noindent{\bf Case 1}. Suppose first that we are in the generic
situation in which $\rho_0({\a})$ is loxodromic and $\tr \rho_0(m)\neq
\pm 2$. Since $\rho_0({\a})$ and $\rho_0(m)$ commute, they have the
same axis. By conjugation we may put the end points of this axis at
$0$ and $\infty$, with the attracting fixed point of $\rho_0({\a})$ at
$\infty$.  For $\rho$ near $\rho_0$, the attracting and repelling
fixed points of $\rho({\a})$ are also holomorphic functions of $\tr
\rho({\a})$, thus we can holomorphically conjugate nearby
representations so that these points are still at $0$ and $\infty$,
respectively. Now choose the local holomorphic branch $f: \R(M) \to
\CC$ of $\l(\rho(m))$ so that $\rho(m) = \pmatrix{e^{f(\rho)} & 0 \cr
  0 & e^{-f(\rho)}} $.  (Notice that the two representations
$\pmatrix{e^{\pm f(\rho)} & 0 \cr 0 & e^{\mp f(\rho)}} $ are conjugate
by rotation by $\pi$ about a point on the axis, but that the
deformations are distinct since there is no smooth family of
conjugations $u_t$ with $u_0 = \mbox{id}$ and $u_t \rho_t^{-1}
u_t^{-1} = \rho_t$.)

It follows easily that under the restriction map
$\pi_m:H^1(M;Ad\rho_0) \rightarrow H^1(m;Ad\rho_0)$, the infinitesimal
deformation $\dot{\rho}$ is mapped to the cocycle $\pi_m(\dot \rho)$
with $\pi_m(\dot \rho)(m)=f'(0)v$, where $f'(0)$ denotes the
derivative of $f(\rho_t)$ at $t=0$.  \medskip

\noindent{\bf Case 2}.
Now suppose that $\rho_0({\a})$ is loxodromic but that $\tr \rho_0(m)
=\pm 2$. Since $\rho_0({\a})$ and $\rho_0(m)$ commute, $\rho_0(m)$
must be the identity in $SL(2,\CC)$.  However, using the fixed points
of $\rho ({\a})$, we can conjugate as before so that $\rho(m)=
\pmatrix{e^{c(\rho)} & 0 \cr 0 & e^{-c(\rho)}} $ for some locally
defined holomorphic function $c$, which we can take to be the {\em
  local definition} of $\l(\rho(m))$.  In fact, since in this
situation both $\rho(\a m)$ and $\rho(\a)$ are loxodromic, one sees
that $c(\rho)$ can be defined by the formula $c(\rho) = \l(\rho(\a m))
- \l(\rho(\a))$.  The discussion then proceeds as before and we again
have $\pi_m(\dot \rho)(m)=c'(0)v$ with $\pi_m$ defined as above.  (In
this case $H^1(m;Id) = \CC^3$.  However, the map $H^1(M; Ad \rho_0)
\to H^1(m; Ad \rho_0)$ factors through $H^1(M; Ad \rho_0) \to
H^1(T_{\a}; Ad \rho_0)$ and one can check directly that $H^1 ( T_{\a};
Ad \rho ) = \CC^2$ is spanned by the cocycles defined by $z_1(m) = v,
z_1({\a}) = 0$ and $z_2(m) =0, z_2(\a) = v$.)

\medskip
\noindent{\bf Case 3}.
Finally suppose that $\rho_0({\a})$ is parabolic.  Since
$\rho_0({\a})$ and $\rho_0(m)$ commute, $\rho_0(m)$ is either
parabolic or the identity. However it is easy to see from the
discussion in the proof of Proposition~\ref{prop:lifting} that, if
$\gamma \in \pi_1(M)$ corresponds to a meridian $m$, then
$\rho_0(\gamma)$ is the product of reflections in the two tangent
circles which contain the plaques of $\dd \cal H $ which meet at the
fixed point of $\rho_0(\gamma)$, and is hence parabolic.

\begin{lemma}
\label{lem:cuspdef} Suppose that $t \mapsto A_t$ is a holomorphic one
parameter family of deformations of a parabolic transformation
$\rho_0(m)$. Then there is a neighborhood $U$ of $0$ in $\CC$ such
that $A_t$ is either
always parabolic or always loxodromic for $t \in U - \{0\}$.
In the first case, $t \mapsto A_t$ is
holomorphically  conjugate to the trivial deformation $t \mapsto A(0)$.
In
the second case,  $A_t$ is holomorphically conjugate to  $
\pmatrix{u_t/2 & (u_t^2 -4)/2 \cr
1/2 & u_t/2} $, where $u_t = \tr A_t$.
\end{lemma}
\begin{proof} Without loss of generality, we may assume that $\tr A_0 =
  2$. For the first statement, note that $t \to \tr A_t$ is
  holomorphic so that $\tr A_t -2$ either vanishes identically or has
  an isolated zero at $0$. Write $A_t = \pmatrix{a(t) & b(t) \cr c(t)
    & d(t)} $.  Conjugating by $z \mapsto 1/z$ if necessary, we may
  assume that $c \neq 0$.  In the first case, translating by $ z
  \mapsto z -(d-a)/c$, we may assume that $A_t= \pmatrix{1 & 0 \cr
    c(t) & 1} $.  In the second case, conjugation by the translation $
  z \mapsto z -(d-a)/2c$ arranges that $a(t)=d(t)$. Conjugating by the
  scaling by $ z \mapsto 4c^2 z$ arranges that $c(t)=1/2$.
\end{proof}
  
It follows that the image of an infinitesimal deformation under the
restriction map is the cocycle which assigns to $m$ the derivative of
$t \mapsto A_t$ with $A_t$ as in the second case in the above lemma.
By direct computation, we calculate that $\pi_m(\dot \rho)(m)=
u'(0)(2u^+ - u^-/4 - v/2)$.  

\smallskip If no elements $\rho_0(m_i)$ are parabolic, then Cases 1
and 2 establish our claim that the restriction map $H^1(M;Ad \rho_0)
\to \oplus_{i=1}^n H^1( m_i;Ad \rho_0)$ is equal to the derivative of
the map $\M$ at $\rho_0$. If some $\rho_0(m_i)$ are parabolic, then
Case 3 shows the same is true provided we replace complex length by
trace.  In summary, we have shown that we can choose, for each $i$, a
linear map $h_i:H^1(m_i;Ad \r) \to \CC$ so that the composition
$(h_1,\ldots, h_n) \circ res:H^1(M; Ad \r) \to \CC^n$ is still
injective and equals the derivative $d{\cal M}_{\r}:H^1(M; Ad \r)
\rightarrow \CC^n$.

\medskip

By similar computations we now show that in the neighborhood of a cusp
the traces of the longitudes can equally be taken as local parameters.
This is crucial in proving Theorem~\ref{thm:main2}.
\begin{prop}
\label{prop:cuspdefok} Let $M$ be a $3$-dimensional
hyperbolic cone manifold and suppose that $\rho_0: \pi_1(M) \to
SL(2,\CC)$ is a lift of the holonomy representation.  For each
boundary torus $T_{\a_j}$, let $m_j$ be the meridian and let $\a_j$ be
the longitude.  Take local parameters $z_j(\rho) = \tr \rho(m_{j})$ if
$\rho_0(m_{j})$ is parabolic and $z_j= \l(\rho(m_{j})) $ otherwise. If
$\rho_0(m_i)$ is parabolic, then $\dd \tr \rho(\a_i) /\dd z_i \neq 0$,
while $\dd \tr \rho(\a_i) /\dd z_j = 0$ for $j\neq i$.
\end{prop}

This result can be extracted from the proof of Thurston's hyperbolic
Dehn surgery theorem, see~\cite{ThuN}.  In fact, $\dd \tr \rho(\a_i)
/\dd \tr \rho(m_{i}) =\tau^2_{i}$, where $\tau_i$ is the modulus of
the induced flat structure on $T_{\a_i}$.  We remark that the Dehn
surgery discussion takes place in a $2^n$-fold covering space of
$\R(M)$ on which one defines complex variables $u_i$ such that
$\tr\rho( m_i) = 2 \ch u_i/2$, see also~\cite{BP} B.1.2.  We shall
give a separate proof which clarifies that
Proposition~\ref{prop:cuspdefok} follows from a fact about
representations of $\pi_1(T^2)$ for a torus $T^2$ into $SL(2,\CC)$. It
is based on the following simple computation:

\begin{lemma}
\label{lem:cuspdef1} Suppose that $t \mapsto \rho_t$ is a holomorphic
one parameter family of deformations of a representation $\rho_0:
\pi_1(T) \to SL(2,\CC)$, defined on a neighborhood $U$ of $0$ in
$\CC$, such that $\rho_0(m) $ is parabolic, $\rho_0(\a) \neq {\rm
  id}$, and such that $\rho_t(m) = A_t$ has the canonical form of
Lemma~\ref{lem:cuspdef} above. Then there exists a holomorphic
function $h: U \to \CC$ such that $h(0) \neq 0$ and such that $
\rho_t(\a) = B_t$ has the form $ \pmatrix{v_t/2 & h(t)(u_t^2 -4)/2 \cr
  h(t)/2 & v_t/2} $, where $v_t = \tr B_t$ and $h^2(t) (u_t^2 -4) =
v_t^2 -4$.
\end{lemma}
\begin{proof} If $t \neq 0$, since $B_t$ and $A_t$ commute, $B_t$ must
  be loxodromic with the same fixed points as $A_t$.  It follows that
  the diagonal entries of $B_t$ must be equal.  Thus $B_t= \pmatrix{
    v_t/2 & k(t)/2 \cr h(t)/2 & v_t/2} $ for analytic functions $h,k$
  with $v_t^2 -4 = hk$. The condition on fixed points gives
  $h^2(t)(u_t^2 -4) =(v_t^2 -4)$ and the form of $B_t$ follows. By
  continuity we must have $B_0= \pmatrix{\pm 1 & 0 \cr h(0)/2 & \pm 1}
  $ and since $\rho_0(\a) $ must be parabolic, $h(0) \neq 0$.
\end{proof}

\noindent{\sc Proof of Proposition~\ref{prop:cuspdefok}:}
To complete the proof, note that the relation $h^2(t) (u_t^2 -4) =
v_t^2 -4$ gives $v'(0) = u'(0) h^2(0)$ which proves the first
statement.  To see that the other derivatives vanish, note that since
$\a_i$ and $m_i$ commute, any deformation which keeps $m_i$ parabolic
necessarily also keeps $\a_i$ parabolic. \ q.e.d.


\section{The local pleating theorem}
\label{sec:localpleating}
In this section we prove Theorem~\ref{thm:localpleating}, the {\em
  local pleating theorem}, which locally characterizes piecewise
geodesic structures by the condition $\tr \sigma(\a) \in \RR$ for all
$\a \in \ua$. This is the first main step in the proof of the local
parameterization Theorem~\ref{thm:main2}.  As usual, let $\bar N$ be a
hyperbolizable $3$-manifold such that $\dd \bar N$ is non-empty and
contains no tori, and let $\ua$ be a doubly incompressible curve
system on $ \dd \bar N$.  We denote by $\cal G( N,\ua)$ the set of
piecewise geodesic structures on $(\bar N,\ua)$ and by $\cal P(N,\ua)$
the subset of convex structures in $\cal G(N,\ua)$. We shall
frequently identify these sets with the corresponding holonomy
representations in $\R(N)$, and topologize $\cal G(N,\ua)$ as a
subspace of $\R(N)$.  Recall that a structure in $ \cal G(N,\ua)$ is
convex if and only if the bending angles satisfy $0 \leq
\theta_{\a}(\sigma) \leq \pi$ for all $\a \in \ua$.

\medskip

We begin with the necessity of the condition that $\sigma(\a)$ have real
trace.
In the case of convex
   structures, this was the starting point of~\cite{KStop}.
\begin{prop}
  \label{prop:realtrace} If $\sigma \in \cal G(N,\u{\a})$, then
$\tr(\sigma(\a)) \in \RR$ for all $\a \in \ua$.
\end{prop}
\begin{proof} This is 
  essentially the same as~\cite{KStop} Lemma 4.6.  Let $(G,\phi)$ be
  the piecewise geodesic structure with holonomy $\phi_*=\sigma$.  If
  $\gamma \in \pi_1(N)$ is freely homotopic to a curve in $\u{\a}$,
  then $g=\s (\gamma)$ is either parabolic, in which case the result
  is obvious, or loxodromic. If $g$ is loxodromic, by definition of a
  piecewise geodesic structure, $\Ax g$ is the intersection of two
  plaques ${\cal N}_1, {\cal N}_2$ of $\phi(\dd \bar N -\ua)$.  Since
  the image of ${\cal N}_i$, $i=1,2$, under $g$ is a plaque which
  contains $\Ax g$, either $g({\cal N}_i)= {\cal N}_i$, $i=1,2$, or
  the two plaques are contained in a common plane $\Pi$ which is
  rotated by $\pi$ and translated along $\Ax g$.  In the first case,
  the half-plane with boundary $\Ax g$ which contains ${\cal N}_1$ is
  mapped to itself under $g$. This can only happen if $g$ is purely
  hyperbolic and hence $\tr \s(\gamma) \in \RR$ as desired. To see
  that the second case cannot arise, consider the $r$-neighborhood $R$
  of $\Ax g$ and its intersection with $\tilde \phi (\bar N - \ua_P)$,
  where $\tilde \phi$ is a lift of $\phi$ to the universal cover of
  $\bar N - \ua_P$.  Since $\tilde \phi$ is an embedding that takes
  $\dd \bar N - \ua_P$ to $\dd \,\mbox{Im}\, \tilde \phi$, we see that
  for small enough $r>0$, $R \cap \mbox{Im}\, \tilde \phi$ is a
  half-tube with boundary $R \cap \Pi$. Since $g$ preserves
  $\mbox{Im}\, \tilde \phi$ and $R$, we see that $g$ cannot rotate
  $\Pi$ by $\pi$.
\end{proof}

In  general, the converse of  Proposition~\ref{prop:realtrace} is
false, see for example Figure 3 in~\cite{KStop}. If however, $\u{\a}$ is
{\em maximal}, the converse holds in the neighborhood of a
convex structure:

\begin{thm}[Local pleating theorem]
\label{thm:localpleating} Let $\sigma_0 \in \cal P(N, \ua)$ where
$\ua$ is a maximal doubly incompressible curve system on $\dd \bar N$.
Let $P$ be the set of elements $\gamma \in \pi_1(N)$ such that
$\sigma_0(\gamma)$ is parabolic.  Then there is a neighborhood $U$ of
$\sigma_0$ in $\R_P(N)$ such that if $\sigma \in U$ and $\tr
\sigma(\a) \in \RR$ for all $\a \in \ua$, then $\sigma \in \cal G(N,
\u{\a})$.
\end{thm}

If the curve system $\ua$ is not maximal, it is easy to see that the
theorem is false, because there are geodesic laminations contained in
what was initially a plaque of $\dd \bar N -\ua$ which are not
contained in $\ua$, along which some nearby structures become bent.
Notice also that although the initial structure is convex, when some
initial bending angle $\theta_\a(\s_0)$ vanishes, we can only conclude
that nearby structures are piecewise geodesic because $\theta_\a(\s)$
can become negative.  If however, the initial bending angles
$\theta_\a(\sigma_0)$ are all strictly positive, the trace conditions
guarantee that locally structures remain convex. A special case of
Theorem~\ref{thm:localpleating} was proved in the context of
quasifuchsian once-punctured tori in~\cite{KSbend}.

\medskip

The idea of the proof is the following.  We always work in a
neighborhood of $[\sigma_0]$ in $\R_P(N)$ in which all groups
$G_{\sigma} = \sigma(\pi_1(N))$ are quasiconformal deformations of
$G_{0} =\sigma_0(\pi_1(N))$. Thus by assumption, $\sigma(\gamma)$ is
parabolic if and only if $\sigma_0(\gamma)$ is parabolic.  Our
assumption that $ \tr \sigma(\a) \in \RR$ implies that if $\gamma \in
\pi_1(N)$ is freely homotopic to a curve in $\ua$, then
$\sigma(\gamma)$ is either parabolic or strictly hyperbolic.  Let
$\cal A_H(\sigma)$ denote the set of axes of the hyperbolic elements
in this set and $\cal A_P(\sigma)$ denote the set of parabolic fixed
points, and let $\cal A (\s) = \cal A_H(\s) \cup \cal A_P(\s)$.

Consider first the group $G_0$. Its convex hull boundary lifts to a
set $X_{\sigma_0} \subset \HH^3$ made up of a union of totally
geodesic plaques which meet only along their boundaries, which are axes
in $\cal A_H(\sigma_0)$. Each component of $X_{\sigma_0}$ separates
$\HH^3$, all the components together cutting out the convex hull $\cal
{H}(G_0)$. (Notice that if $\dd \bar N$ is compressible,
$X_{\sigma_0}$ may not be simply connected.  Nevertheless, the closure
of exactly one component of $\HH^3 - X_{\sigma_0}$ contains $\cal
H(G_0)$.)

Now suppose we have $\s$ near $\s_0$ such that $\tr \s(\a) \in \RR$
for $\a \in \ua$.  The axes $\cal A_H(\sigma)$ are near to those in
$\cal A_H(\sigma_0)$. Because the traces remain real, axes in a common
plaque remain coplanar, so that we can define a corresponding union of
plaques $X_{\sigma}$. The main point is to show that, like the plaques
making up $X_{\sigma_0}$, these nearby plaques also intersect only
along their boundaries, in the corresponding axes of $\cal A_H(\s)$.
In other words, with the obvious provisos about smoothness along the
bending lines, $X_{\sigma}$ is a $2$-manifold without boundary
embedded in $\HH^3$. Then a standard argument can be used to show that
each component of $X_{\sigma}$ separates $\HH^3$. Together the
components cut out a region $E_{\s}$ which is close to the convex hull
$\cal H(G_0)$. Finally we show that the quotient $E_{\s}/G_{\s}$ is
the image of the induced embedding $\phi_{\s}:\bar N - \ua_P
\hookrightarrow N(G_{\s})$. This defines a piecewise geodesic
structure on $(\bar N, \ua)$ with parabolic locus $\ua_P$. 

\medskip

In more detail we proceed as follows.  First consider the initial
convex structure $\phi_0:\bar N -\ua_P \hookrightarrow N(G_0)$ with
holonomy representation $\sigma_0$.  Since $\ua$ is maximal, the
closure of each component $Q$ of $\phi_0(\dd \bar N-\ua)$ is a totally
geodesic pair of pants with geodesic boundary (where we allow that
some of the boundary curves may be punctures), so a lift $\tilde
Q=\tilde Q_{\sigma_0}$ of such a component $Q_{\sigma_0}$ will be
contained in a plane $\Pi(\tilde Q) \subset \HH^3$. Let $\Gamma(\tilde
Q)$ be the stabilizer of $\Pi(\tilde Q)$ in $G_0$. The closure of
$\tilde Q$ in $\Pi(\tilde Q)$ is the Nielsen region (i.e. the convex
core) $\cal N(\tilde Q)$ of $\Gamma(\tilde Q)$ acting on $\Pi(\tilde
Q)$; by definition $\cal N(\tilde Q)$ is a plaque of $\phi_0(\bar N -
\ua_P)$. Since $ Q$ is a three holed sphere (where a hole may be a
puncture), $\Gamma(\tilde Q)$ is generated by three suitably chosen
elements $\sigma_0(\gamma_i), i=1,2,3$ whose axes project to the three
boundary curves of the closure of $Q$.

Let $X_{\sigma_0}$ be the union of all the Nielsen regions.  Since
$\phi_0$ is a convex structure, $\phi_0 (\dd\bar N - \ua_P) = \dd
V(G_0)$ and so $ X_{\sigma_0} = \dd {\cal H}(G_0)$.  Each plaque $\cal
N(\tilde Q)$ is adjacent to another plaque $\cal N(\tilde Q')$ along
an axis in $\cal A_H(\sigma_0)$. Moreover, if $\cal N(\tilde Q), \cal
N(\tilde Q')$ are distinct plaques then their intersection is either
empty or coincides with an axis in $\cal A_H(\sigma_0)$. Thus $
X_{\sigma_0}$ is a $2$-manifold without boundary in $\HH^3$.

Now suppose we have a representation $[ \sigma ]$ near $[ \sigma_0 ]$
in $\R_P(N)$. By normalizing suitably, we can arrange that
$\sigma(\gamma)$ is arbitrarily near $\sigma_0(\gamma)$ for any finite
set of elements $\gamma \in \pi_1(N)$.  The assumption is that $\tr
\sigma(\gamma) \in \RR$ whenever $\gamma$ is freely homotopic to a
curve in $\u{\a}$. This implies that if $ \sigma_0(\gamma_1),
\sigma_0(\gamma_2), \sigma_0(\gamma_3)$ generate $\Gamma(\tilde
Q_{\sigma_0})$, then the subgroup $\Gamma(\tilde Q_\sigma)$ generated
by $\sigma(\gamma_1), \sigma (\gamma_2),\sigma (\gamma_3)$ is Fuchsian
with invariant plane $\Pi(\tilde Q_\sigma)$ (see for
example~\cite{Indra} Project 6.6).  Here $\tilde Q_\sigma$ is the
interior of the Nielsen region $\cal N(\tilde Q_\sigma)$ of
$\Gamma(\tilde Q_\sigma)$ acting on $\Pi(\tilde Q_\sigma)$. Define
$X_\sigma$ to be the union of all the Nielsen regions $\cal N(\tilde
Q_\sigma)$.  Without presupposing that the structure $\s$ is piecewise
geodesic, call $ \cal N(\tilde Q_\sigma) $ a {\em plaque} of
$X_\sigma$. Note that $ \cal N(\tilde Q_\sigma) $ is determined by the
axes $\Ax \sigma(\gamma_i)$ if $\sigma(\gamma_i)$ is hyperbolic, or
the fixed points and tangent directions of $\sigma(\gamma_i)$ if
$\gamma_i$ is parabolic where $i=1,2,3$.

As sketched above, we want to show that the regions making up
$X_{\sigma}$ intersect only along their boundaries, in other words,
that $X_{\sigma}$ is a $2$-manifold embedded in $\HH^3$.  We begin
with a lemma which describes how distinct plaques can intersect.  Let
$Hom_P(\pi_1(N),SL(2,\CC))$ denote the subset of $\s \in Hom(\pi_1(N),
SL(2,\CC))$ such that $\s(\gamma)$ is parabolic for all $\gamma \in
P$. For convenience, we denote by $\RR(\ua)$ the subset of elements
$\s \in Hom(\pi_1(N),SL(2,\CC))$ satisfying the condition that $\tr
\s(\a) \in \RR$ for all $\a \in \ua$.
\begin{lemma}
\label{lemma:intersection}  
For $\s$ near $\s_0$ in $Hom_P(\pi_1(N),SL(2,\CC)) \cap \RR(\ua)$, if
two distinct plaques $\cal N(\tilde Q_\sigma) $,$\cal N(\tilde
Q'_\sigma)$ intersect, then the intersection either coincides with an
axis in $\cal A_H(\sigma)$ or must meet such an axis.
\end{lemma}
\begin{proof}
  Suppose first that the two plaques intersect transversely.  Since
  each plaque is planar, their intersection is a geodesic arc
  $\hat{\beta}$ which either continues infinitely in at least one
  direction, ending at a limit point in $\Lambda(\Gamma(\tilde
  Q_\sigma)) \cap \Lambda(\Gamma(\tilde Q'_\sigma))$, or which has
  both endpoints on axes in $\cal A_H(\s)$. In the first case, since
  $\Lambda(\Gamma(\tilde Q_\sigma)) \cap \Lambda(\Gamma(\tilde
  Q'_\sigma))=\Lambda(\Gamma(\tilde Q_\sigma) \cap \Gamma(\tilde
  Q'_\sigma)) \neq \emptyset$ (see for example~\cite{MT} Theorem
  3.14), we have that $\Upsilon_{\sigma}=\Gamma(\tilde Q_\sigma) \cap
  \Gamma(\tilde Q'_\sigma) \neq \{1\}$.  Since $\Upsilon_{\sigma}$
  preserves both $\cal N(\tilde Q_\sigma)$ and $\cal N(\tilde
  Q'_\sigma)$, it preserves the geodesic segment $\hat \beta$ in which
  they intersect.  Therefore, $\Upsilon_{\sigma}$ is an elementary
  subgroup generated by a hyperbolic isometry whose axis $\beta$
  contains $\hat \beta$.
  
  Now, for $\s$ near $\s_0$, since $ \sigma\sigma_0^{-1}:G_{0}
  \rightarrow G_\sigma$ is a type-preserving isomorphism which maps
  $\Gamma(\tilde Q_{\sigma_0})$, $\Gamma(\tilde Q'_{\sigma_0})$ to
  $\Gamma(\tilde Q_\sigma)$, $\Gamma(\tilde Q'_\sigma)$ respectively,
  it follows that $\Upsilon_0=\Gamma(\tilde Q_{\sigma_0}) \cap
  \Gamma(\tilde Q'_{\sigma_0})$ is also generated by a loxodromic
  isometry.  Its axis must lie in both of the Nielsen regions $\cal
  N(\tilde Q_{\sigma_0}) $ and $\cal N(\tilde Q'_{\sigma_0})$ and must
  therefore be a geodesic in $\cal A_H(\sigma_0)$. Thus, in this case,
  $\hat \beta$ must continue infinitely in both directions so that
  $\hat \beta=\beta$ and $\beta$ must be contained in $\cal
  A_H(\sigma)$.

Finally, if $\cal N(\tilde Q_\sigma) $ and $\cal N(\tilde Q'_\sigma)$
are coplanar, the same argument works if we choose $\hat\beta$ to be
{\em any} geodesic in $\cal N(\tilde Q_\sigma) \cap \cal N(\tilde
Q'_\sigma)$.
\end{proof}

The point of the above lemma is that intersections between plaques
always meet in the inverse image of a suitably chosen compact subset
of $X_{\s}/G_\s$, because we can always arrange for the axes $\cal
A_H(\sigma)$ not to penetrate far into the cusps. More precisely, by
the Margulis lemma, for each $\sigma$ we can choose a set of disjoint
horoball neighborhoods of the cusps in $\HH^3/G_\sigma$. If $p_\sigma$
is a parabolic fixed point of $G_\sigma$, let $H(p_\sigma)$ denote the
corresponding lifted horoball in $\HH^3$. Since we are deforming
through type preserving representations, we may assume that in a
neighborhood $U$ of $\sigma_0$, the horoballs $H(p_\sigma)$ vary
continuously with $\sigma$, in the sense that in the unit ball model
of $\HH^3$, their radii and tangent points move continuously.
Moreover, since the finitely many geodesics whose lifts constitute
$\cal A_H(\sigma)$ have uniformly bounded length in $U$, they
penetrate only a finite distance into any cusp.  Therefore, by
shrinking the horoballs $H(p_\sigma)$ and replacing $U$ by a smaller
neighborhood if necessary, we may assume that $\cal A_H(\sigma) \cap
H(p_\sigma) = \emptyset$ for all $p_\sigma$ and for all $\sigma \in
U$.  Thus, the lemma implies that if two plaques intersect, then their
intersection meets in $Y_\s =X_\sigma \cap H_\s$, where $H_\s= \HH^3 -
\cup H(p_\s)$ is the complement of the horoball neighborhoods.

 \medskip

 The action of $\Gamma(\tilde Q_{\sigma_0})$ on $\cal N(\tilde
 Q_{\sigma_0})$ has a fundamental polygon (for example, made of two
 adjacent right angled hexagons, where some of the sides may be
 degenerate if $\Gamma(\tilde Q_{\sigma_0})$ contains parabolics)
 whose intersection $F(\tilde Q_{\sigma_0})$ with $H_{\sigma_0}$ is
 compact. Choose a fundamental polygon for each pair of pants and let
 $K_0=\cup_{i=l}^k F_i(\sigma_0)$ be the union of such compact pieces,
 where $k$ is the total number of pairs of pants in $\dd \bar N -
 \ua$.  We can define for $\s$ near $\s_0$, corresponding fundamental
 polygons and compact set $K_\s=\cup_{i=l}^k F_i(\sigma)$.  The
 projection $K_\s/G_\s$ is equal to $Y_\s/G_\s$.  We shall denote the
 plaque containing $F_i(\s)$ as ${\cal N}_i(\s)$. Clearly, the
 projection $\cup_{i=1}^k {\cal N}_i(\s)/G_\s$ is equal to $X_\s
 /G_\s$. In particular, an arbitrary plaque ${\cal N}(\s)$ of $\s$ is
 a translate $\s(\gamma) {\cal N}_i(\s)$ for some $i \in \{1,\ldots,k\}$
 and some $\gamma \in \pi_1(N)$.
 
 Let $\epsilon>0$ and define $K$ to be the closed
 $\epsilon$-neighborhood of $K_0$ in $\HH^3$.  Since each compact set
 $F_i(\s)$ is determined by a finite number of axes and parabolic
 points in $\cal A(\s)$ and since the position of an axis or parabolic
 point in $\cal A(\sigma)$ varies continuously with $\sigma$, there
 exists a neighborhood $U$ of $\s_0$ such that the Hausdorff
 distance between $F_i(\sigma)$ and $F_i(\sigma_0)$ is at most
 $\epsilon$ for all $\sigma \in U$ and for all $i=1,\ldots,k$.
 Thus, $K_\s \subset K$ for all $\s \in U$.
\begin{prop}
\label{prop:plaques}
There exists a neighborhood $U$ of $\sigma_0$ in
$Hom_P(\pi_1(N),SL(2,\CC))$ with the property that for $\sigma \in U \cap
\RR(\ua)$, two distinct plaques ${\cal N}(\sigma),{\cal
  N}'(\sigma)$ intersect only along axes in ${\cal A}_H(\sigma)$.
Thus $X_\s$ is a $2$-manifold without boundary in $\HH^3$.
\end{prop}
\begin{proof}
   Suppose there were no such neighborhood. Then there exists a
   sequence of representations $\s_n \rightarrow \s_0$ and pairs of
   plaques ${\cal N}(\s_n),{\cal N}'(\s_n)$ which intersect along
   geodesic segments which are not contained in axes in ${\cal
     A}_H(\sigma_n)$. By Lemma~\ref{lemma:intersection}, ${\cal N}(\s_n)
   \cap {\cal N}'(\s_n)$ has non-empty intersection with $H_{\s_n}$.
   By translating if necessary, we may therefore assume that $ {\cal
     N}(\s_n) \cap {\cal N}'(\s_n) \cap K_{\s_n} \neq \emptyset$. By
   taking a subsequence of $\s_n$ if necessary, we can further assume
   that ${\cal N}(\s_n) \cap {\cal N}'(\s_n) \cap F_i(\s_n) \neq
   \emptyset$ for some $i$ and for all $n$.

   Let us consider ${\cal N}(\s_n)$. It will become clear that the
   following line of argument can also be applied to ${\cal N}'(\s_n)$.
   Since ${\cal N}(\s_n)$ is a translate of one of ${\cal N}_1(\s_n),
   \ldots,{\cal N}_k(\s_n)$, we can take a further subsequence of
   $\s_n$ if necessary and assume that ${\cal N}(\s_n)$ is a translate
   of ${\cal N}_j(\s_n)$ for some $j$, for all $n$. In other words,
   there exists a sequence of elements $\gamma_n \in \pi_1(N)$ such
   that ${\cal N}(\s_n)=\s_n(\gamma_n) {\cal N}_j(\s_n)$. It follows
   that $\s_n(\gamma_n) {\cal N}_j(\s_n) \cap F_i(\s_n) \neq
   \emptyset$.  Furthermore, by composing $\gamma_n$ with another
   deck-transformation if necessary and using the fact that $\s_n$
   preserves $H_{\s_n}$, we can assume that
\begin{equation}
\label{eqn:plaques}
\s_n(\gamma_n) F_j(\s_n) \cap F_i({\s_n})  \neq \emptyset.
\end{equation}

Now let $K$ be the compact set defined in the discussion preceding the
statement of the proposition. Since $\s_n \rightarrow \s_0$, we have
that $F_j(\s_n), F_i(\s_n) \subset K_{\s_n} \subset K$ for large $n$.
Then Equation(\ref{eqn:plaques}) automatically implies that
$$\s_n(\gamma_n) K \cap K \neq \emptyset.$$

Since the set $\{g \in SL(2,\CC): g(K) \cap K \neq \emptyset \}$ is
compact, by passing to a subsequence, we may assume that
$\sigma_n(\gamma_n) \rightarrow g_0$ for some $g_0 \in SL(2,\CC)$.
Thus $g_0$ is contained in the geometric limit of the groups
$\sigma_n(\pi_1(N))$. However, since $\sigma_0(\pi_1(N))$ is
geometrically finite and since $\s_n$ is type preserving, the
convergence is strong, see for example~\cite{MT} Theorem 7.39
or~\cite{Kap} Theorem 8.67. Thus $\sigma_n(\gamma_n) \rightarrow
\sigma_0(\gamma)$ for some $\gamma \in \pi_1(N)$, and hence
$\sigma_n(\gamma_n \gamma^{-1}) \rightarrow \mbox{id}$. Since
$\sigma(\pi_1(N))$ is always discrete, we have
$f(\sigma)=\inf\{d(\sigma(\delta),\mbox{id}): \delta \in \pi_1(N),
\sigma(\delta) \neq \mbox{id} \}>0$. Therefore, by choosing a small
enough neighborhood $U$ of $\sigma_0$, we can guarantee that $f$
restricted to $U$ is bounded below by a strictly positive constant. It
then follows that $\sigma_n(\gamma_n \gamma^{-1})=\mbox{id}$ for large
$n$, in other words, $$\sigma_n(\gamma_n) = \sigma_n(\gamma). $$

Since $F_j(\s_n), F_i(\s_n)$ converge to $F_j(\s_0), F_i(\s_0)$
respectively, the preceding, together with
Equation(\ref{eqn:plaques}), imply that $\s_0(\gamma)F_j(\s_0) \cap
F_i(\s_0) \neq \emptyset$ and so $\s_0(\gamma){\cal N}_j(\s_0) \cap
{\cal N}_i(\s_0) \neq \emptyset$.  Now for $\s_0$, we know that any
two plaques which intersect either coincide or intersect in an axis in
${\cal A}_H(\s_0)$.  Therefore, either $\s_0(\gamma){\cal N}_j(\s_0)=
{\cal N}_i(\s_0)$ or $\s_0(\gamma){\cal N}_j(\s_0) \cap {\cal
  N}_i(\s_0)$ is an axis in ${\cal A}_H(\s_0)$. The first case implies
that ${\cal N}(\s_n)= {\cal N}_i(\s_n)$, for large $n$. The second
case implies that ${\cal N}(\s_n)\cap {\cal N}_i(\s_n)$ is an axis in
$A_H(\s_n)$, for large $n$.

Since the same argument can be also applied to ${\cal N}'(\s_n)$, by
comparing to the intersections for $\s_0$, we deduce that for large
$n$, ${\cal N}(\s_n)$ and ${\cal N}'(\s_n)$ either coincide or
intersect along an axis in $A_H(\s_n)$, both of which contradict the
hypothesis. 
\end{proof}

\begin{cor}
\label{cor:separates} Each component of $X_\sigma$ separates $\HH^3$.
\end{cor}
\begin{proof} This is a standard topological argument,
   see for example~\cite{Hirsch} Theorem 4.6.  Let $X$ be a component
   of $X_\s$. Note that if $x_1,x_2 \notin X$ then the mod $2$
   intersection number $I(\beta)$ of a path $\beta$ joining $x_1$ to
   $x_2$ with $X$ is a homotopy invariant, which moreover only depends
   on the components of $\HH^3 - X$ containing $x_1$ and $x_2$. Since
   $\HH^3$ is simply connected, $I$ is constant.  If $X$ did not
   separate, $I$ would be even. However by choosing points close to
   opposite sides of a plaque of $X$, we see $I$ is odd.
\end{proof}

\noindent {\sc Proof of Theorem~\ref{thm:localpleating}.}
For the convex structure $\sigma_0$, choose a point $x_0 \in \HH^3$ in
the interior of the convex hull $\cal H(G_0)$ which projects to the
thick part of $\HH^3/G_0$. Since $X_{\sigma_0}= \dd \cal H(\sigma_0)$
and since $X_\sigma$ moves continuously with $\sigma$, we may assume
that $x_0 \notin X_\sigma$ for $\sigma$ near $\sigma_0$.  For each
component $X^i = X^i_\sigma$ of $X_\sigma$, let $E^i = E^i_\sigma$ be
the closure of the component of $\HH^3 - X_i$ which contains $x_0$. If
$X^j$ is another component of $X_\sigma$, then $X^j \subset
\mbox{Int}E^i$ and we argue as in Corollary~\ref{cor:separates} that
$X^j$ separates $E^i$. Hence, for $\sigma$ near $\sigma_0$, the set
$E_{\sigma} = \cap_i E^i_\sigma$ is non-empty. By construction
$E_{\sigma}$ is $G_{\sigma}$-invariant and closed.  Notice also that
no end of $\HH^3/G_0$ is contained in $E_{\sigma}/G_{\sigma}$.

We claim that $E_{\sigma}/G_{\sigma}$ is homeomorphic to $\bar N -
\ua_P$.  Suppose first that $\theta_{\a_i}(\sigma) \geq 0$ for all
$i$. In this case we actually have equality
$E_\sigma/G_\sigma=V(G_\sigma)$. To see this, first note that
$E^i_\sigma$ is locally convex and therefore it is convex
(see~\cite{CEG} Corollary 1.3.7). Thus $E_{\sigma}$ contains the
convex hull $\cal H(G_{\sigma})$.  Moreover, by construction
$X_\sigma$ is contained in the convex span of $\cal A_H(\sigma) \cup
\cal A_P(\sigma)$, so that $X_\sigma = \dd E_{\sigma} \subset \cal
H(G_\sigma)$. If $ \cal H(G_{\sigma}) \neq E_{\sigma}$, then there is
a point $x \in \dd \cal H(G_{\sigma}) \cap \mbox{Int} E_{\sigma}$.
Since we are assuming that $G_{\sigma}$ is geometrically finite, there
is a bijective correspondence between components of $ \dd \cal
H(G_{\sigma})$ and lifts of ends of $\HH^3/G_{\s}$, which in turn
correspond to components of the regular set $\Omega(G_{\s})$.  Thus if
$x$ is in a component $Z_{\s}$ of $ \dd \cal H(G_{\sigma})$, we can
find a geodesic arc $\gamma$ starting from $x$ and ending on $\dd
\HH^3$ in the component $\Omega^i_{\s}$ of $\Omega(G_{\s})$ which
`faces' $Z_{\s}$, and such that points on $\gamma$ near $x$ are not in
$ \cal H(G_{\sigma})$.

Since the corresponding $X^i_{\s_0}$ separates $Z_{\s_0}$ from
$\Omega^i_{\s_0}$, it follows that $X^i_{\s}$ separates $Z_{\s}$ from
$\Omega^i_{\s}$. Thus $\gamma$ must intersect $X^i_{\s}$. Since
$X^i_{\s} \subset \cal H(G_{\sigma})$, this gives a geodesic subarc of
$\gamma$ with endpoints in $ \cal H(G_{\sigma})$ parts of whose
interior are outside $ \cal H(G_{\sigma})$, contradicting convexity.
It follows that $E_{\sigma}= \cal H(G_{\sigma})$ and hence in this
case we have $E_\sigma/G_\sigma=V(G_\sigma)$.

Now we consider the general case where $\theta_{\a_i}(\s)$ may be
negative for some $i$. Use the ball model of $\HH^3$ and let
$B^3=\HH^3 \cup \dd \HH^3$. First, observe that the closure
$\overline{\cal H(G_{0})}$ of $\cal H(G_0)$ in $B^3$ is a closed ball
whose boundary is the union of $\dd \cal H(G_0)$ and the limit set
$\Lambda_{\s_0}$ of $G_0$. Next, consider the closure $\overline E_\s$
of $E_\s$ in $B^3$.  We shall prove below that $\dd \overline E_\s=
\dd E_\s \cup \Lambda_\s$, where $\dd E_\s $ is as usual the boundary
in $\HH^3$.  Assuming this fact, let us show that $\dd \overline E_\s$
is an embedded $2$-sphere in $B^3$.  On the one hand, each component
of $\dd E_\s$ is homeomorphic to a corresponding component of $\dd
\cal H(G_0)$ by a homeomorphism $h_{\s}$ which varies continuously
with $\s$. On the other hand, the $\lambda$-lemma~\cite{MSS} gives the
analogous result for the limit set $\Lambda_\s$. More precisely, there
is an open neighborhood $W$ of $\s_0$ in $Hom_P(\pi_1(N),SL(2,\CC))$
and a continuous map $f:\Lambda_{\s_0} \times W \rightarrow \dd B^3$
such that $f(\xi,\sigma_0)=\xi$ and that $f(\,\cdot\,,\s)$ is a
homeomorphism $\Lambda_{\s_0} \rightarrow \Lambda_\s$ for all $\s \in
W$.  Using equivariance, it is easy to check that these homeomorphisms
glue together to induce a homeomorphism between $\dd \overline E_\s=
\dd E_\s \cup \Lambda_\s$ and $\dd \overline{\cal H(G_0)} = \dd \cal
H(G_0) \cup \Lambda_{\s_0}$.  We deduce that $\dd \overline E_\s$ is
an embedded $2$-sphere in $B^3$ as claimed.

Since $B^3$ is irreducible, $\overline E_\s$ must be a $3$-ball. Thus
$E_\s$ is the universal cover of $E_\s/G_\s$. Since $\pi_1(E_\s/G_\s)
\approx G_\s \approx \pi_1(\bar N -\ua_P)$ and since by construction
$\dd E_\s$ projects to a union of surfaces homeomorphic to $\dd \bar N
- \ua_P$, we can apply Waldhausen's Theorem~\cite{Wald} to conclude
that there is a homeomorphism $\phi_\s: \bar N -\ua_P \rightarrow
E_\s/G_\s$ which induces $\s$ as required.

Finally, we prove our claim that $\dd \overline E_\s= \dd E_\s \cup
\Lambda_\s$, or equivalently that $\dd \overline E_\s \cap \dd B^3 =
\Lambda_\s$.  First, note that $\dd \overline E_\s \cap \dd B^3$ is
closed and $G_\s$-invariant in $\dd B^3$ and therefore contains the
limit set $\Lambda_\s$ of $G_\s$.  If $z \in \dd \overline E_\s \cap
\dd B^3$, pick $z_n \in E_\s$ with $z_n \to z$.  Let $\Delta$ be a
closed fundamental domain for $E_\s$.  Choose $g_n \in G_\s$ with $w_n
= g_nz_n \in \Delta$.  If $w_n$ is eventually in $ \Delta \cap Y_{\s}$
then by compactness we may assume that $w_n \to w \in \Delta \cap
Y_{\s}$. Then the sequence $ g_n^{-1}w$ also accumulates on $z$, so
that $z \in \Lambda_\s$.

Otherwise, $w_n$ and hence $z_n$ is eventually contained in the union
$\cal U$ of the horoball neighborhoods $H(p_{\s})$. Note that the
intersection of $\cal H(G_{\s})$ with $H(p_{\s})$ is the region bounded
between two planes tangent at $p_{\s}$, and $\overline {\cal H(G_{\s})}
\,\cap \, \overline {H(p_{\s})} \cap \dd B^3 = \{p_{\s}\}$. Since no axes
in $\cal A_H(\s)$ meet $H(p_{\s})$, from the definition we have $E_{\s}
\cap H(p_{\s}) = {\cal H}(G_{\s}) \cap  H(p_{\s}) $. Hence $\dd
\overline{E_{\s} }\,\cap \, \overline {H(p_{\s})} \cap \dd B^3 =
\{p_{\s}\}$.  We deduce that the limit of any sequence eventually
contained in $\cal U$ is either a parabolic point or a limit of
parabolic points, and is therefore in $\Lambda_{\s}$.  This completes
the proof that $\dd \overline E_\s \cap \dd B^3 = \Lambda_\s$ as
required.  \ q.e.d. \vskip .2in

\noindent{\bf Remark} It is easy to see that in fact, in general,
$E_\s \subset \cal H(G_\s)$. From the above, $\overline E_\s$ is a
ball whose boundary $\dd \overline E_\s= \dd E_\s \cup \Lambda_\s$ is
contained in $\overline{\cal H(G_\s)}=\cal H(G_\s) \cup \Lambda_\s$.
Since $\overline{\cal H(G_\s)}$ is itself a closed $3$-ball, it
follows that the interior of $\overline E_\s$ is entirely contained in
$\overline{\cal H(G_\s)}$ and hence that $E_\s \subset \cal H(G_\s)$.


\section{Local isomorphism of representation spaces}
\label{sec:isomorphism}

In this section we prove that if $\ua$ is a {\em maximal} doubly
incompressible curve system on $\dd \bar N$, then $\cal R(M)$ and
$\cal R(N)$ are locally isomorphic near a convex structure on $(\bar
N,\ua)$.  Here, as usual, $M=D\bar N - \ua$ as in
Section~\ref{sec:doubles}.  This is the second main step in the proof
of the local parameterization theorem, Theorem~\ref{thm:main2}.

\begin{thm}[Local isomorphism theorem]
\label{thm:isomorphism} Let $\ua=\{\a_1,\ldots,\a_d\}$ be a maximal
curve system on $\dd \bar N$.  Let $\sigma_0 \in \R(N)$ be a convex
structure in $\cal P(N,\u{\a})$ and let $\rho_0 \in \R(M)$ be its
double. Then the restriction map $r: \R(M) \rightarrow \R(N)$,
$r(\rho) = \rho|_{\pi_1( N )}$, is a local isomorphism in a
neighborhood of $\rho_0$.
\end{thm}
By Theorems~\ref{thm:smoothreps} and~\ref{thm:kh}, $\R(N)$ and $\R(M)$
are both complex manifolds of dimension $d$ at $\sigma_0 $ and
$\rho_0$ respectively. Thus it will suffice to prove that $r$ is
injective. To do this, we also consider the natural restriction map $
{\hat r}: \R(M) \rightarrow \R({\hat N}), \ \hat r(\rho) =
\rho|_{\pi_1(\tau(N))}$, where $\tau$ is the involution on $M$ and
$\R({\hat N})$ is the representation variety of the manifold $\hat N =
\tau(N)$. We will prove the injectivity of $r$ by factoring through
the product map $(r,{\hat r}): \R(M) \rightarrow \R(N) \times \R({\hat
  N})$.  Thus we first consider the effect of a deformation of
$\rho_0$ on the induced structures on both halves $N$ and ${\hat N}$,
and then show that the symmetry of $M$ implies that what happens on
${\hat N}$ is fully determined by what happens on $N$.
\begin{prop}
\label{prop:injectivity} Let $\ua$ be a maximal doubly incompressible
curve system on $\dd \bar N$. Suppose that $\sigma_0 \in \P(N,\ua )$,
and let $\rho_0 \in \R(M)$ be its double. Then the restriction map
$(r,{\hat r}): \R(M ) \rightarrow \R(N ) \times \R({\hat N} )$ is
injective on a neighborhood of $\rho_0$.
\end{prop}
\begin{proof}
  It will suffice to prove that the derivative of $(r,{\hat r})$ is
  injective on tangent spaces.  As explained in
  Section~\ref{subsec:complexlength}, we can identify the tangent
  spaces to $\R(N)$ and $\R(M)$ with the cohomology groups
  $H^1(\pi_1(N);\Ad \sigma_0)$ and $H^1(\pi_1( M);\Ad \rho_0)$
  respectively. Thus showing that the induced map $$(r_*,{\hat r}_*):
  T_{\r}\R(M) \rightarrow T_{r(\r)}\R(N) \times T_{{\hat
      r}(\r)}\R({\hat N})$$
  is injective is the same as showing the
  induced map on cohomology is injective. We claim that it will be
  sufficient to show that if a cocycle $z \in Z^1(\pi_1(M);Ad \rho_0)$
  satisfies the condition that $z(\gamma)=0$ and $z(\hat \gamma)=0$
  for all $\gamma \in \pi_1(N)$, $\hat \gamma \in \pi_1(\hat N)$, then
  $z \equiv 0$. To see why this is so, first note that if $z$ induces
  the $0$-class in $H^1(\pi_1(N);Ad \,r(\rho_0))$, we can modify $z$
  by a coboundary so that $z(\gamma)=0$ for all $\gamma \in \pi_1(N)$.
  Now, since we are assuming that $z$ also induces the $0$-class in
  $H^1(\pi_1(\hat N);Ad \,\hat r(\rho_0))$, we have that $z(\hat
  \gamma)=v - Ad \rho_0(\hat \gamma)v$ for all $\hat \gamma \in
  \pi_1(\hat N)$, where $v$ is some fixed element in $\sl2c$. 
  We need to see that $v = 0$.
  
  Since
  $\ua$ is maximal, all components of $\dd \bar N - \ua$ are pairs of
  pants; we call the loops round their boundaries pants curves. All
  loops under consideration will have a fixed base point $x_0$, which
  we choose so that it lies in one of the pants $Q_0$ in a component
  $S_0$ of $\partial \bar N $.  We will use $\gamma$ to denote both a
  loop and its representative in $\pi_1(N,x_0)$. Note that a loop
  $\gamma$ completely contained in $N $ has a mirror loop
  $\tau_*(\gamma)$ in ${\hat N} $, as $\tau$ interchanges $N $ and
  ${\hat N} $ and fixes $\partial \bar N - \ua= \partial {\overline
    {\hat N}} - \ua$ pointwise.  Let $\gamma_1,\gamma_2$ be two pants
  curves of $Q_0$.  Since the involution $\tau$ fixes both
  $\gamma_1,\gamma_2$, we have that
  $$0=z(\gamma_i)=z(\tau_*(\gamma_i))=v-Ad \rho_0(\tau_*(\gamma_i))v=
  v-Ad \rho_0(\gamma_i)v,$$
  for $i=1,2$. However, since the two
  isometries $\rho_0(\gamma_1),\rho_0(\gamma_2)$ do not commute, it
  must be that $v=0$.
  
  Suppose then that $z(\gamma)=z(\hat \gamma)=0$ for all $\gamma \in
  \pi_1(N), \hat \gamma \in \pi_1(\hat N)$.  Since $\r$ is a cone
  structure on $M$, by Theorem~\ref{thm:kh}, any infinitesimal
  deformation will be detected by an infinitesimal change in the
  holonomy of some meridian curve.  Thus to show that $z \equiv 0$, it
  will be sufficient to show that $z(m_i)=0$ for every meridian $m_i$,
  $i=1,\ldots,d$.

  We will first show that the deformations induced on the meridians
  associated to the pants curves of $Q_0$ are trivial. Choose homotopy
  classes $m_{\a}$ of the meridians as depicted in
  Figure~\ref{fig:curves}.  There may be two different types,
  depending on whether or not the corresponding boundary curve is
  shared by a different pair of pants.
\begin{figure}[htb]
\centerline{\epsfbox{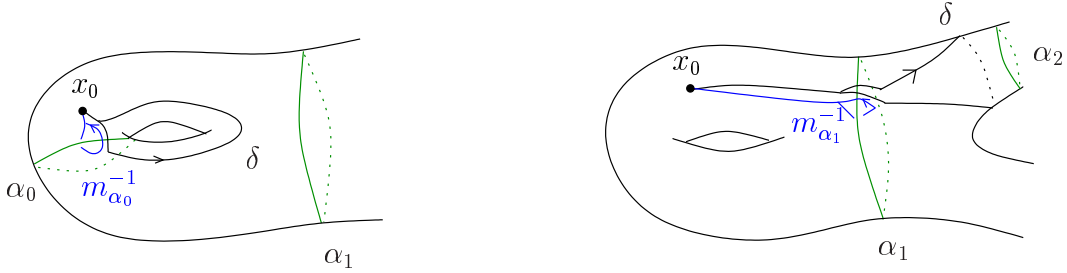}}
\caption{Homotopy classes of meridians.}
\label{fig:curves}
\end{figure}

If $\alpha=\alpha_0$ is a pants curve which is not shared by another
pair of pants, as on the left in Figure~\ref{fig:curves}, then it has
a dual curve $\delta$ which intersects it only once and intersects
none of the other curves in $\ua$.  Hence:
$$\tau_*(\delta)=m_{\a}^{-1} \cdot \delta.$$
Then for any cocycle $z \in Z^1(M; Ad \rho_0)$, we have
$$
z(\mera  \cdot \tau_*(\delta))=z(\mera) + Ad\r(\mera)z(\tau_*(\delta))=
z(\delta).
$$
Since $z(\delta)$ and $z(\tau_*(\delta))$ are both zero by
assumption, it must be that $z(\mera)$ is also zero.

The other possibility for a pants curve $\a$ of $Q_0$ is that it is
shared by an adjacent pair of pants, such as $\a_1$ on the right in
Figure~\ref{fig:curves}. Such an $\a$ has a dual curve $\delta$ which
intersects it exactly twice and intersects none of the other curves in
$\ua$. In particular, we can choose $\delta$ to be freely homotopic to
a pants curve in the adjacent pants. Hence, we have the relation
$$ \tau_* (\delta)=\mera^{-1} \cdot \delta \cdot \mera$$
from which we obtain
\begin{eqnarray}
&& z(\delta \cdot \mera)= z(\delta) + Ad\r(\delta)z(\mera) \nonumber
\\
&=& z(\mera \cdot \tau_*(\delta)) = z(\mera)+
Ad\r(\mera)z(\tau_*(\delta)). \nonumber
\end{eqnarray}
Since $z(\delta)$ and $z(\tau_*(\delta))$ are both zero by hypothesis,
$z(\mera)$ must be contained in the centralizer of $\r(\delta)$. On
the other hand, since $m_{\alpha}$ and $\a$ commute, $z(m_{\a})$ is
also in the centralizer of $\r(\a)$. However, since $\r(\a)$ and
$\r(\delta)$ do not commute, $z(m_{\a})$ must be zero.

We now proceed by an inductive argument on adjacent pairs of pants.
For this purpose, we choose homotopy classes for the meridians in a
tree-like fashion, as shown in Figure~\ref{fig:tree}, first focusing
on the component $S_0$.
\begin{figure}[htb]
\centerline{\epsfbox{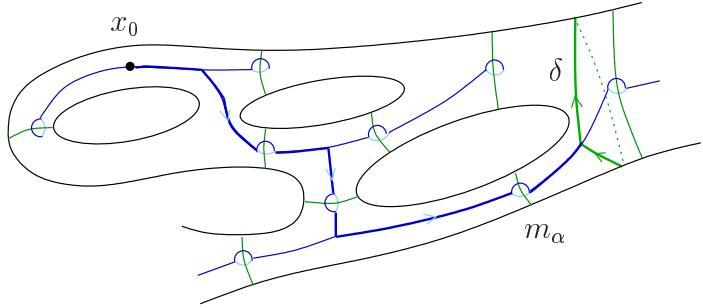}}
\caption{A tree of meridians.}
\label{fig:tree}
\end{figure}

Let $\a \in \ua$ be the boundary of some pair of pants $Q \subset
S_0$. Our inductive hypothesis is that there is a chain of pants
$Q_0,Q_1,\ldots, Q_j=Q$ contained in $S_0$ such that each meridian
associated to a pants curve in $Q_i$, $i<j$, has trivial deformation.
To show that $z(m_{\a})=0$, we wish to choose a curve $\delta$ dual to
$\a$ and apply the same argument as before. However, since all loops
are based at $x_0$, any dual curve is forced to intersect a collection
of pants curves. We choose $\delta$ so that it meets a succession of
pants curves contained in $\cup_{i \leq j} Q_i$, say
$\a_1,\a_2,\ldots,\a_k$ as indicated in the figure. Then
\begin{equation}
\label{eqn:relation}
\tau_*(\delta) = m_{\a_1}^{-1} \cdots m_{\a_k}^{-1} \cdot \delta \cdot
m_{\a_k} \cdots m_{\a_1}.
\end{equation}
We have that $z(m_{\a_1}), \ldots, z(m_{\a_{k-1}})$ are zero by the
inductive hypothesis and that $z(\tau_*(\delta)), z(\delta)$ are zero
by the underlying assumption.  Denote the product $m_{\a_{k-1}} \cdots
m_{\a_1}$ by $x$. Then Equation(\ref{eqn:relation}) gives $x \cdot
\tau_*(\delta) \cdot x^{-1} =m_{\a_k}^{-1} \cdot \delta \cdot
m_{\a_k}$ which implies that $z(m_{\a_k}^{-1} \cdot \delta \cdot
m_{\a_k}) = 0$, and hence $ z(m_{\a_k})=Ad\r(\delta)z(m_{\a_k}).$ For
the same reason as before, $z(m_{\a_k})=z(m_\alpha)$ is zero.  Thus we
have shown that $z(m_\alpha)=0$ for all $\alpha$ in $S_0$.

The same method can be used to show that $z(m_{\a})$ is zero for all
$\a$ in any other component $S$ of $\partial \bar N $. Choose a pair
of pants $Q$ in $S$ and fix a point $x_1$ in $Q$. Let $s$ be an arc in
$N $ from $x_1$ to $x_0$. If $\b$ is a loop based at $x_1$, then the
concatenation $s * \b * s^{-1}$ is a loop based at $x_0$. We can now
repeat the previous arguments using loops of this form.
Observe that for each pants curve $\a_j$ in $Q$, the loop
$\kappa= s* \tau(s)^{-1}$ satisfies the relation
$$\tau_*(\a_j)=\kappa^{-1} \cdot \a_j \cdot \kappa$$
This implies that $z(\kappa)=0$. Using this fact, we can show that
$z(m_{\a_j})=0$ for the pants curves $\a_j$ in $Q$ and then apply the
inductive argument. In place of Equation(\ref{eqn:relation}) we have
relations of the form
$$\tau_*(\delta)=
\kappa^{-1} \cdot m_{\a_1}^{-1} \cdots m_{\a_k}^{-1} \cdot \delta \cdot
m_{\a_k} \cdots m_{\a_1} \cdot \kappa.$$
However, since $z(\kappa)=0$, the calculations are identical.
\end{proof}

\medskip

The following proposition, which exploits the symmetry between $N$ and
$\hat N$, now completes the proof of Theorem~\ref{thm:isomorphism}.
Denote by $\cal C(M,\ua)$ the space of cone structures on $M$ with
singularities along $\ua$.  We also write $\mu_{\a}$ for the complex
length $\l(\rho(m_{\a}))$ of the meridian $m_{\a}$.  The crucial
observation is that since $\cal C(M,\ua)$ is the purely imaginary
locus of the coordinate functions $\mu_{\a}$, a holomorphic function
on $\cal R(M)$ is locally determined by its values on $\cal C(M,\ua)$.

\begin{prop}
\label{prop:injectivity2}
Let $\sigma_0$ be a convex structure in $\cal P( N, \u{\a})$ and let
$\rho_0$ be its double. Then in a neighborhood of $(r(\rho_0),\hat
r(\rho_0))$, the projection $ \R(N ) \times \R({\hat N} ) \rightarrow
\R(N ) $ is injective on the image $(r,\hat r )(\R(M ))$.
\end{prop}
\begin{proof} Lifting $\rho_0$ to an element in  $ \mbox{Hom}(M,
  SL(2,\CC))$, from the construction in
  Proposition~\ref{prop:lifting}, we have $\rho_0 \circ \tau_* =
  J\rho_0J^{-1}$, where $J$ is a reflection in a plane in $\HH^3$. By
  considering first the reflection induced by $J_0(z)=\bar z$, it is
  easy to check that $\tr JAJ^{-1} = \tr \overline{A}$ for any
  orientation reversing isometry $J$ of $\HH^3$ and $A \in SL(2,\CC)$.
  Thus $\tr \rho_0 \circ \tau_* (\gamma) = \tr\overline{\rho_0
    (\gamma)}$ for all $\gamma \in \pi_1(M)$, where $\overline{A}$ is
  the matrix whose entries are complex conjugates of those of $A$. Let
  $\overline{\rho_0}$ be the representation defined by
  $\overline{\rho_0}(\gamma)=\overline{\rho_0(\gamma)}$.  This shows
  that the two representations $\rho_0 \circ \tau_*$ and
  $\overline{\rho_0}$ are conjugate and thus are equivalent in
  $\R(M)$, see for example\,\cite{CuS}.  The main point of the proof
  is to show that for all cone structures $\rho$ near $\rho_0$, we
  have
\begin{equation}
\label{eqn:equalityoncone}
\rho \circ \tau_*= \overline \rho
\end{equation}
as equivalence classes in $\R(M)$, where $\overline{\rho}(\gamma)=
\overline{\rho(\gamma)}$.

First consider the generic case where $\rho_0(\a)$ are loxodromic
for all $\a \in \ua$. Then, by Theorem~\ref{thm:kh},
a holomorphic deformation of $\rho_0$
is parameterized by the complex lengths $\mu_j$ of the meridians
$m_j$, $j=1,\ldots,d$. We shall show below that
\begin{equation}
\label{eqn:symmetry}
\mu_j(\rho \circ \tau_*)=\mu_j(\bar \rho),
\end{equation}
for all $j=1,\ldots,d$ and $\rho \in \C(M,\ua)$ near $\rho_0$,
which implies Equation(\ref{eqn:equalityoncone}).

Let $\a, m \in \pi_1(M)$ be commuting representatives of a longitude
and meridian pair $\a_j,m_j$. Following the discussion in
Section~\ref{subsec:complexlength}, in order to compute the complex
length $\mu_{\a}$ of $ \rho(m)$, we first conjugate $\rho$ so that the
axis of the longitude $\rho({\a})$ is in {\em standard position},
meaning that its repelling and attracting fixed points are at
$0,\infty$, respectively.  The matrix $A= \rho(m)$ is then diagonal
and $\mu_{\a} (\rho)$ is the logarithm of the top left entry.  Now
consider the representation $\bar \rho$. Notice that when $\rho$ is
conjugated so that $\rho({\a})$ is in standard position, so is ${\bar
  \rho}(\a)$.  We can therefore read off the complex length
$\mu_{\a}({\bar \rho})$ of ${\bar \rho} (m)$ from the matrix $\bar A=
\overline{\rho(m)}$ (see Cases 1 and 2 in
Section~\ref{subsec:complexlength}).  We deduce that $\mu_j(\bar
\rho)=\overline{\mu_j(\rho)}$ for all $j=1,\ldots,d$.  On the other
hand, since $\tau_*(m)$ is conjugate to $m^{-1}$ in $\pi_1(M)$ by the
same element which conjugates $\tau_*(\a)$ to $\a$, we can use the
same method to also deduce that $\mu_j(\rho\circ \tau_*)=-\mu_j
(\rho).$ Now for $\rho \in \C(M,\ua)$, the functions $\mu_j(\rho)$ are
purely imaginary so that $\overline{\mu_j(\rho)}=-\mu_j(\rho)$.  The
three equalities give Equation(\ref{eqn:symmetry}) as desired.

We must also consider the case in which some of the meridians are
parabolic. For these meridians, the parameter in question is their
trace. Note that
$$\tr\rho(m_j)=\tr\rho(m_j^{-1}) =
\tr \rho \circ \tau_*(m_j),$$
and
$$\tr \bar \rho(m_j)=\overline {\tr \rho(m_j)} .$$
In particular, $
\tr \rho \circ \tau_*(m_j)= {\tr \bar \rho(m_j)} $ whenever $ \tr \rho
(m_j) \in \RR$. (Notice we do not need to
assume that all manifolds $ \tr \rho (m_j) \in \RR$ are cone
manifolds; in fact the meridian $\rho (m_j)$ will be purely hyperbolic
for some points in the real trace locus near the parabolic point.)

In summary, if we take local coordinates $w_j(\rho) =\tr \rho
(m_j)$ whenever $\rho_0(m_j)$ is parabolic and $w_j(\rho) =
\sqrt{-1}\mu_j (\rho)$ otherwise, we have shown in the two cases above
that for all $j=1,\ldots,d$
\begin{equation}
w_j(\rho \circ \tau_*) = w_j(\overline \rho)
\end{equation}
on the $d$-dimensional real submanifold of $\R(M)$ locally defined by
the condition $(w_1,\ldots,w_d) \in \RR^d$. Since the map $\rho
\mapsto \rho \circ \tau_*$ is holomorphic in $\rho$, this implies that
for a holomorphic deformation $\rho_t$ of $\rho_0$, where $t$ is a
complex variable in a neighborhood of $0$, we have
\begin{equation}
\label{eqn:equalityall}
[\rho_t \circ \tau_*]
=[\overline{\rho_{\overline{t}}}]
\end{equation}
as equivalence classes in $\R(M)$. (Since we are concerned with 
deformations only up to first order, we are assuming here that $t$ is
real if and only if $w_j(t)$ is real, for all $j$.)

We now wish to equate the cocycles defined by the two deformations.
Recall that we are identifying the cohomology group $H^1(\pi_1(M);Ad
\,\rho)$ with the {\em holomorphic} tangent space to $\R(M)$ at
$\rho$. In particular, this means that
$$\rho_t(\gamma)=\left( Id + t \dot{\rho}(\gamma) + O(|t|^2) \right)
\rho_0(\gamma)$$
where $\dot \rho =z \in Z^1(\pi_1(M);Ad \,\rho_0)$ is the cocycle
defined by
$$z(\gamma)=\
\frac{d}{dt}\Big|_{t=0}\rho_t(\gamma)\rho_0(\gamma)^{-1}.$$
Therefore,
$$\overline{\rho_{{\overline t}}}(\gamma)= \left(Id  + t
   \overline{\dot{\rho}(\gamma)} + O(|t|^2) \right)
\overline{\rho_0}(\gamma).$$
In other words, the cocycle $w$
associated to $\overline{\rho_{\overline{t}}}$ has values given by
$w(\gamma)=\overline{z(\gamma)}$.  We emphasize that if $z$ is a
cocycle in $Z^1(\pi_1(M); Ad\,{\rho_0})$, then the function $\bar z$
whose values are given by $\bar z(\gamma) = \overline{z(\gamma)}$, is
naturally a cocycle in $Z^1(\pi_1(M); Ad\,\overline{\rho_0})$.

On the other hand, the cocycle associated to $\rho_t \circ \tau_*$ is
given by
$$\frac{d}{dt}\Big|_{t=0}\rho_t(\tau_*(\gamma))\rho_0(\tau_*(\gamma))^{-1}=
z(\tau_*(\gamma)).$$
Again, note that $z \circ \tau_*$ is naturally a
cocycle in $Z^1(\pi_1(M);Ad\, \overline{\rho_0})$.

Thus, it follows from Equation(\ref{eqn:equalityall}) that
${z \circ \tau_*}$ and $\overline{z}$ differ by a coboundary
in $B^1(\pi_1(M);Ad\, \overline{\rho_0})$. In
other words, for all $\gamma \in \pi_1(M)$, we have
$${z(\tau_*(\gamma))}= \overline{z(\gamma)} + v - Ad\,
\overline{\rho_0}(\gamma)v$$
where $v$ is some element in $\sl2c$.
Hence, if $z(\gamma)=0$ for all $\gamma \in \pi_1(N)$, it follows that
$$z(\tau_*(\gamma))={v} - Ad \overline{\rho_0(\gamma)}
{v}={v} - Ad {\rho_0(\tau_*(\gamma))}{v},$$
i.e. $z(\hat \gamma)={v} - Ad {\rho_0(\hat \gamma)}{v}$
for all $\hat \gamma \in \pi_1(\hat N)$.
\end{proof}
This concludes the proof of the
local isomorphism theorem~\ref{thm:isomorphism}.

We single out the following useful fact extracted from the above
proof, see also~\cite{Otal} Section 3.
\begin{cor}
\label{cor:realtrace}
Let $\s_0$ be a convex structure in $\P(N,\ua)$ and let $\rho_0$ be
its double. Then there is a neighborhood $U$ of $\rho_0$ in $\cal
C(M,\ua)$ such that $\rho \circ \tau_*=\bar\rho$ for all $\rho \in U$.
In particular, $\tr\rho(\gamma)=\overline{\tr \rho(\gamma)}$ whenever
$\gamma$ is freely homotopic to a curve on $\dd \bar N$.
\end{cor}

\section{Lengths are parameters}
\label{sec:main}
\subsection{Local parameterization of $\R(N)$}

We begin by proving the {\em local} parameterization theorem,
Theorem~\ref{thm:main2}.  To make a precise statement, we first
clarify the definition of the complex length map $\cal L:\R(N)
\rightarrow \CC^d$.  Let $\ua=\{\a_1,\ldots,\a_d\}$ be a maximal
doubly incompressible curve system on $\dd \bar N$ and let $\s_0 \in
\P(N,\ua)$ be a convex structure.  Number the curves in $\ua$ so that
$\s_0(\a_i)$ is parabolic for $i=1,\ldots,k$ and purely hyperbolic
otherwise.  Define
$$\L(\sigma)=(\tr \s(\a_1),\ldots, \tr
\s(\a_k),\l_{\a_{k+1}}(\s),\ldots,\l_{\a_d}(\s)),$$
where $\l_{\a_{i}}(\s)$ denotes the complex length $\l(\s(\a_i))$.
We can then state Theorem C as:
\begin{thm}
\label{thm:main1}
Let $\ua$ be a maximal curve system on $\dd \bar N$ and let $\s_0 \in
\P(N,\ua)$ be a convex structure.  Then $\L:\R(N) \rightarrow \CC^d$
is a local holomorphic bijection in a neighborhood of $\sigma_0$.
  \end{thm}

We will actually show:
\begin{thm}
\label{thm:main1'} The composition
${\cal L} \circ r: \R(M) \rightarrow \CC^d$ is a local holomorphic
bijection in a neighborhood of the double $\rho_0$ of $\s_0$, where
$r:\R(M) \rightarrow \R(N)$ is the restriction map.
\end{thm}

Combined with Theorem~\ref{thm:isomorphism}, this second result gives 
Theorem~\ref{thm:main1}.

\medskip 

Recall that $\s_0(\a_i)$ is parabolic if and only if $\rho_0(m_i)$ is
parabolic for the associated meridian $m_i$.  (As usual, we use $\s
(\a_i)$ to mean $\s (\gamma)$ for $\gamma \in \pi_1(N)$ freely
homotopic to $\a_i$.)  Theorem~\ref{thm:kh} implies that
$(z_1,\dots,z_d)$ are local coordinates for $\R(M)$ near $\rho_0$,
where $z_i = \tr \rho(m_i)$ for $i \leq k$ and $z_i=\l_{m_i}(\rho)=
\l(\rho(m_i))$ for $i>k$.  Split the Jacobian of ${\cal L}\circ r$ at
$\rho_0$ into four blocks by cutting the matrix between the rows
$k,k+1$ and between the columns $k,k+1$. Let $m=d-k$.
Proposition~\ref{prop:cuspdefok} says that the lower off-diagonal
block of size $m\times k$ is the $0$-matrix and that the $k\times k$
block $\left(\dd \tr r(\rho)(\a_i) )/ \dd \tr \rho(m_j) \right)$ is a
diagonal matrix none of whose diagonal entries vanish.  Therefore, to
show that the Jacobian is non-singular, it is sufficient to show that
the $m \times m$ block ${\cal J}=\left( \dd \l_{\a_i}(r(\rho))/ \dd
  \l_{m_j}(\rho) \right)$ is non-singular.

Observe that ${\cal J}$ is the Jacobian of the map
$F:\R_P(M)\rightarrow \CC^m$ defined by $$F(\rho)=
(\l_{\a_{k+1}}(r(\rho)),\ldots,\l_{\a_d}(r(\rho))),$$
where $\R_P(M)$
is the set representations $\rho$ for which $\rho(m_i)$ is parabolic
for $i \leq k$.  (By Theorem~\ref{thm:kh}, $\R_P(M)$ is a smooth
complex manifold of dimension $m$ near $\rho_0$ parameterized by the
complex lengths $(\l_{m_{k+1}}(\rho),\ldots,\l_{m_d}(\rho))$.)
Clearly $F$ is the composition ${\cal L}_m \circ r_P$ of the
restriction map $r_P:\R_P(M) \rightarrow \R_P(N)$ and ${\cal L}_m:
\R_P(N) \rightarrow \CC^m$, where ${\cal
  L}_m(\s)=(\l_{\a_{k+1}}(\s),\ldots,\l_{\a_d}(\s))$.

The crucial observation is that $F$ is a `real map' with respect to
the two totally real $m$-submanifolds $\cal C \cap \R_P(M)$ in the
domain and $\RR^m$ in the range, where ${\cal C}={\cal C}(M,\ua)$
denotes the cone structures in $\R(M)$ with singular locus $\ua$.
This is where both the local pleating theorem~\ref{thm:localpleating}
and the local isomorphism theorem~\ref{thm:isomorphism} are used:

\begin{prop}
\label{prop:forwardbackwardopen}
There is an open
neighborhood $V$ of $\rho_0$ in $\cal R_P(M)$ such that $F(V \cap {\cal
C (M,\ua)})$
is contained in $ \RR^m$; and there is an open
neighborhood $U$ of $F(\r)$ in $\CC^m$ such that $F^{-1}(U \cap \RR^m)$
is contained in ${\cal C} (M,\ua)$.
\end{prop}
\begin{proof} The first statement is immediate from
  Corollary~\ref{cor:realtrace}.  By the local pleating
  theorem~\ref{thm:localpleating}, there is a neighborhood $W$ of
  $F(\rho_0)$ in $\CC^m$ such that ${\cal L}^{-1}(\RR^m \cap W)$ is
  contained in $\cal G(N,\ua)$.  Now for each structure $\sigma \in
  {\cal L}^{-1}(\RR^m \cap W)$ near $\sigma_0$, there is a cone
  structure $\rho \in {\cal C}(M,\u{\a})$ near $\r$ such that
  $r(\rho)=\sigma$, namely, its double. One deduces easily from
  Theorem~\ref{thm:isomorphism} that the restriction $r:{\cal R} (M)
  \rightarrow {\cal R} (N)$ induces a local isomorphism $r:{\cal
    R}_P(M) \rightarrow {\cal R}_P(N)$ in a neighborhood of $\r$. It
  follows that we can find a neighborhood $U$ of $\sigma_0$ in ${\cal
    R}_P(N)$ such that $r^{-1}(\cal G(N,\u{\a}) \cap U) \subset {\cal
    C} (M,\u{\a})$. The result follows.\end{proof}

We complete the proof of Theorem~\ref{thm:main1'} using a result from
complex analysis.  If $H: \CC \to \CC$ is a holomorphic map in one
variable such that $H(0) =0$, $H(\RR) \subset \RR$ and $H^{-1}(\RR)
\subset \RR$, then it is easy to see that $H$ is non-singular at $0$.
The following shows that this result extends to holomorphic maps of
$\CC^m$.

\begin{prop}
\label{prop:holo}
Let $Z,W $ be open neighborhoods of ${\bf 0} \in \CC^m$ and suppose
that $H: Z \to W$ is a holomorphic map such that $H({\bf 0}) ={\bf
  0}$, $H(Z \cap \RR^m ) \subset \RR^m$ and $H^{-1}(W \cap \RR^m)
\subset \RR^m$. Then $H$ is invertible on a neighborhood of ${\bf 0}$.
\end{prop}
{\bf Proof.} We will prove that $\mbox{d}H$ is non-singular at ${\bf
  0}$.  The key is that $H$ cannot be branched and therefore must be
one-to-one.

Take coordinates $(z_i)$ for $Z$ and $(w_i)$ for $W$.  First consider
the complex variety $H^{-1}\{{\bf 0}\}$. By hypothesis, it is
contained in $\RR^m$. This is not possible unless $H^{-1}\{{\bf 0}\}$
is a variety of dimension $0$, i.e., a discrete subset of points, for
otherwise, the coordinate functions $(z_i)$ would be real-valued on
the complex variety $H^{-1}\{{\bf 0}\}$. In this case, $H$ is said to
be {\em light} at ${\bf 0}$ and there are open neighborhoods ${\tilde
  U}$ of ${\bf 0}$ in the domain and $U$ of ${\bf 0}$ in the range
such that the restriction $H|_{\tilde U}: {\tilde U} \rightarrow U$ is
a finite map (see~\cite{Loj}, Section V.2.1).  Furthermore, by
Remmert's Open Mapping Theorem, $H$ is an open map.

If $\mbox{d}H({\bf 0})$ is singular, then the matrix $\left(
  \frac{\partial w_i}{\partial z_j}({\bf 0})\right)$ is singular. Let
$u_i = \mbox{Re}\, w_i$ and $x_i = \mbox{Re}\, z_i$.  Since $H(\tilde
U \cap \RR^m) \subset \RR^m$ it follows that $\frac{\partial
  u_i}{\partial x_j}({\bf 0}) \in \RR$. Moreover $\frac{\partial
  w_i}{\partial z_j}({\bf 0})= \frac{\partial w_i}{\partial x_j}({\bf
  0}) = \frac{\partial u_i}{\partial x_j}({\bf 0})$, and hence the
real matrix $A=\left( \frac{\partial u_i}{\partial x_j}({\bf
    0})\right)$ is singular.  In particular, there is a real line
$\ell$ in $\RR^m$ whose tangent vector is not contained in the image
of $A$. Let $D$ denote the corresponding complex line in $\CC^m$. Then
$V=H^{-1}(D) \cap {\tilde U} $ is a $1$-dimensional complex variety.

Now, using the classical local description of real and complex $1$
dimensional varieties as in~\cite{Mil} Lemma 3.3, we can pick a branch
$ {\tilde D}$ of $V$ which is locally holomorphic to $\CC$ and such
that ${\tilde D} \cap \RR^m$ has a single branch locally holomorphic
to $\RR$.  Then $h=H|_{\tilde D}: {\tilde D} \rightarrow D$ is a
non-constant holomorphic map from one complex line to another. Since
the image of $\mbox{d}h({\bf 0})$ is trivial, $h$ must be a branched
covering of degree $k$ with $k>1$. However, by hypothesis,
$h^{-1}(\ell)$ is contained in $\RR^m$ and hence ${\tilde D} \cap
\RR^m$ contains a union of $k$ distinct lines, a contradiction.  
\ q.e.d. \vskip .2 in

\noindent{\sc Proof of Theorem~\ref{thm:main1'}.}
   Proposition~\ref{prop:forwardbackwardopen} shows that, after
translating origins, the  map
$F $ satisfies the conditions
on $H$ in Proposition~\ref{prop:holo}. Since $r$ is a bijection by
Theorem~\ref{thm:isomorphism}, the result follows.
\ q.e.d. \vskip .2in
This completes the proofs of Theorems~\ref{thm:main1}
and~\ref{thm:main1'}.

\subsection{Global parameterization of $\cal P(N,\ua)$}
In this last section, we prove the global parameterization
Theorems A,B stated in the introduction.
\begin{prop}
\label{prop:limit}
Let $\ua$ be a maximal doubly incompressible curve system. Then
there is a unique point $\sigma_{*}$ in $\cal P(N,\ua)$ such that
$\theta_{i} = \pi$ for all $i$. For any $\sigma_0 \in \cal P(N,\ua)$,
there is path $\sigma_t \in \cal P(N,\ua),\ t \in [0,1]$ with initial
structure $\sigma_0$ such that $\sigma_1 = \sigma_{*}$.
\end{prop}
\begin{proof}
  Since $\ua$ is maximal, it follows from 
  Proposition~\ref{prop:sufnonemptypleating} and~\cite{KMS} that there is a
  unique hyperbolic structure $\s_*$ on $N$ in which all curves in
  $\ua$ are parabolic and thus where all bending angles are $\pi$.

  Let $\sigma_0 \in \P(N,\ua)$ have bending angles
  $(\theta_{1}{(\sigma_0)},\ldots,\theta_{n}(\sigma_0))$.  Suppose
  first that $\sigma_0 \in \P^+(N,\ua)$, so that
  $\theta_{i}{(\sigma_0)}>0$ for all $i$.  By
  Theorem~\ref{thm:boglobalpar} the bending angles parameterize
  $\P^+(N,\ua)$ and form a convex set. Hence there is a $1$-parameter
  family of structures $\s_t \in \P^+(N,\ua)$ defined by $
  \theta_{i}(\sigma_t)= \theta_{i}(\sigma_0) + t(
  \pi-\theta_{i}(\sigma_0))$ for each $i$, where $0 \leq t \leq 1$.
  This clearly defines the required path.
   
   Now assume that some of the initial angles vanish. The fact that
   $\P^+(N,\ua)$ is locally parameterized by the bending angles is
   deduced from the Hodgson-Kerckhoff theorem (our
   Theorem~\ref{thm:kh}) in Lemme 23 of~\cite{BonO}. We need the
   extension of this result to $\P(N,\ua)$.  As long as the cone
   manifold obtained by doubling $N$ has non-zero volume, the
   Hodgson-Kerckhoff theorem allows that some of the cone angles may
   be $2\pi$, equivalently that some of the bending angles may vanish.
   Based on this observation, an inspection of the proof of Lemme 23
   shows that the local parameterization for $\P^+(N,\ua)$ goes
   through unchanged to $\P(N,\ua)$. In other words, the bending
   angles are local parameters in a neighborhood of $\s_0$.
   Therefore, regardless of whether some initial bending angles are
   zero, there is a small interval $[0,\epsilon]$ for which the
   structures $\s_t$, $t \in [0,\epsilon]$ are uniquely determined by
   $\s_0$ and are contained in $\P(N,\ua)$. Since
   $\theta_i(\s_\epsilon)>0$ for all $i$, we are now in the situation
   in which Theorem~\ref{thm:boglobalpar} applies and we proceed as
   before.
\end{proof}
\begin{prop}
\label{prop:posdef} Suppose that $\ua$ is maximal and that $\sigma_0
\in \cal P(N,\ua)$. Let
   $\varphi_{\a_j}=2(\pi-\theta_{j})$ be the cone-angle along $\a_j$.
Then the Jacobian matrix
$${\mbox d}L(\sigma_0)=\left( \frac{\partial l_{\a_i}}{\partial
\varphi_{\a_j}}
(\sigma_0)\right)$$
is positive definite and symmetric.
\end{prop}
\begin{proof} Renumber the curves $\a$ so that $\sigma_0(\a_i)$ is
  parabolic for $i=1,\ldots, k$ and purely hyperbolic otherwise and
  let $\rho_0$ be the double of $\s_0$. For nearby $\rho$, take local
  coordinates $(z_1,\ldots, z_d)$ where $z_i=\tr \rho(m_i)$ for $i \le
  k$ and $z_i=\l_{m_i}(\rho)$ otherwise. Likewise, near $\s_0$, take
  local coordinates $(w_1,\ldots,w_d)$ where $w_i=\tr \s(\a_i)$ for $i
  \le k$ and $w_i=\lambda_{\a_i}(\s)$ otherwise. Although for $i\leq
  k$, the complex length $\lambda_{\a_i}$ cannot be defined in a
  neighborhood of $\rho_0$, we can always pick a branch of
  $\lambda_{\a_i}$ so that on $\cal P(N,\ua)$ it is a non-negative
  real valued function $l_{\a_i}$, coinciding with the hyperbolic
  length of $\a_i$.  In Theorem~\ref{thm:main1} we showed that the
  matrix $\mbox{d}\L(\s_0)=\left( \frac{\dd w_i}{\dd z_j}(\s_0)
  \right)$ is non-singular. To show that $\mbox{d} L(\s_0)$ is
  non-singular, we compare the entries of the two matrices.
    
  In Proposition~\ref{prop:cuspdefok} we showed that the upper left $k
  \times k$ submatrix of $\mbox{d}\L(\s_0)$ is diagonal with non-zero
  entries.  In fact, the diagonal entries are strictly negative. We see
  this as follows.  From Lemma~\ref{lem:cuspdef1}
  $$h_i(\rho)^2(\tr \rho(m_i)^2 -4)= (\tr \rho(\a_i)^2 -4)$$
  for a
  locally defined holomorphic function $h_i(\rho)$ with the property
  that $h_i(\rho_0) \neq 0$.  For cone structures $\rho$ near
  $\rho_0$, it follows from Corollary~\ref{cor:realtrace} that both
  $\tr \rho(m_i)^2 -4$ and $\tr \rho(\a_i)^2 -4$ are real-valued.  A
  careful inspection shows that $\tr \rho(m_i)^2 -4$ and $\tr
  \rho(\a_i)^2 -4$ must be of opposite sign, for otherwise, in the
  limit, $h_i(\rho_0)$ would be real, making the holonomies
  $\rho_0(m_i)$ and $\rho_0(\a_i)$ parabolic with the same translation
  direction. However this contradicts the fact that the translation
  directions are orthogonal, as discussed in the beginning of Case 3
  in Section~\ref{subsec:complexlength}. Thus for $i \leq k$,
  $$\frac{\dd w_i}{\dd z_i}(\rho_0) =\frac{\dd \tr \rho(\a_i)}{\dd \tr
    \rho(m_i)}\Big|_{\rho=\rho_0} =h_i^2(\rho_0) <0.$$
  Since $$\frac{\dd w_i}{\dd z_i}(\rho_0)= - \frac{\partial
      l_{\a_i}}{\partial \varphi_{\a_i}}(\sigma_0),$$
  it follows that the upper left $k\times k$ submatrix of
  $\mbox{d}L(\sigma_0)$ is diagonal with strictly positive entries.

Since deformations which keep $\rho(m_i)$ parabolic also keep
$\rho(\a_i)$ parabolic, for $i \leq k$ and $j>k$ we have
$$\frac{\dd w_i}{\dd z_j}(\sigma_0)=0=\frac{\partial
  l_{\a_i}}{\partial \varphi_{\a_j}}(\sigma_0).$$
For $j \leq k$ and
$i>k$ we calculate directly that for points $\sigma \in \P(N,\ua)$
near $\sigma_0$ we have $z_j=2\cos \frac{\varphi_{\a_j}}{2}$ so
    $$
    \frac{\partial l_{\a_i}}{\partial \varphi_{\a_j}} =
    \frac{\partial l_{\a_i}}{\partial z_j} \frac{\partial
      z_j}{\partial \varphi_{\a_j}} = -\frac{\partial
      l_{\a_i}}{\partial z_j} \sin \frac{\varphi_{\a_i}}{2} \to 0.$$

Finally, for $i,j>k$, we have $\varphi_{\a_j}=\mbox{Im}\,z_j$ and
$l_{\a_i}=\mbox{Re}\,w_i$. Recall that if $\rho \in \C(M,\ua)$, then
$w_i=l_{\a_i}$ at $r(\rho)$. Therefore, for $i,j >k$,
$$ \frac{\dd w_i}{\dd z_j}(\sigma_0)= -\sqrt{-1}
\frac{\partial l_{\a_i}}{\partial \varphi_{\a_j}}(\sigma_0).$$
It follows that the matrix $\mbox{d}L(\sigma_0)$ is also non-singular.

Now observe that the Schl\"{a}fli formula for the
volume of the convex core~\cite{Bon} gives
$$\mbox{dVol}=-\frac{1}{2}\sum_{i=1}^{n}l_{\a_i}\mbox{d}\varphi_{\a_i}.$$
Thus the matrix $\mbox{d}L(\sigma_0)$ is the Hessian of the volume
function on $\P(N,\ua)$. This automatically implies the symmetry
relation
$$ \frac{\partial l_{\a_i}}{\partial \varphi_{\a_j}}=
    \frac{\partial l_{\a_j}}{\partial \varphi_{\a_i}},$$
for all $i,j>k$.

This discussion shows that when $\sigma=\sigma_*$ represents the
structure in which all bending angles are $\pi$, the matrix
$\mbox{d}L(\sigma_*)$ is diagonal and that all diagonal entries are
strictly positive.  In particular, it is positive definite and
symmetric. By Proposition~\ref{prop:limit}, $\sigma_*$ can be
connected to any given $\sigma_0 \in \P(N,\ua)$ by a path $\sigma_t$
in $\P(N,\ua)$. Since, by the same reasoning as above,
$\mbox{d}L(\sigma_t) $ is non-degenerate along this path, it must
remain positive definite, proving our claim.
\end{proof}
\begin{cor}
\label{cor:posdef}
Let $\ua = \{\a_1,\ldots,\a_n\}$ be any doubly incompressible curve
system on $\bar N$, not necessarily maximal. Let $\s_0 \in \P(N,\ua)$
and let $\varphi_{\a_j}=2(\pi - \theta_{\a_j})$ be the cone-angle
along $\a_j$. Then the matrix
$$\left(\frac{\partial l_{\a_i}}{\partial \varphi_{\a_j}}(\sigma_0)
\right)_{i,j \leq n}$$ is positive definite and symmetric.
\end{cor}
\begin{proof}
  Extend $\ua$ to a maximal system $\ua' $ by adding curves
  $\{\a_{n+1},\ldots,\a_d\}$. By Proposition~\ref{prop:posdef}, the 
enlarged matrix
$$\left(\frac{\partial l_{\a_i}}{\partial \varphi_{\a_j}}(\sigma_0)
\right)_{i,j \leq d}$$
is positive definite and symmetric.  Since a
symmetric submatrix of a positive definite symmetric matrix is itself
positive definite, the claim follows.
\end{proof}
\noindent {\bf Theorem B (Mixed parameterization)}
{\em For any ordering of the curves $\ua$ and for any $q$, the map $\sigma
\mapsto (l_{\alpha_1}(\s), \ldots, l_{\alpha_q}(\s),
\theta_{\alpha_{q+1}}(\s), \ldots,\theta_{\alpha_{n}}(\s))$ is 
an injective local diffeomorphism on $\P(N,\ua)$.}

\medskip

\noindent \begin{proof}
  Corollary~\ref{cor:posdef} shows that the map is a local
  diffeomorphism, so we have only to show it is injective.  Suppose
  there are two points $\s_0,\s_1 \in \P(N,\ua)$ such that
  $l_{\a_i}(\s_0)=l_{\a_i}(\s_1)$ for $i \leq q$ and
  $\theta_{\a_i}(\s_0)=\theta_{\a_i}(\s_1)$ for $i>q$.  To simplify
  notation, let $v_i = \theta_{\a_i}(\s_0)$ and $u_i =
  \theta_{\a_i}(\s_1)$ for all $i=1,\ldots,n$. It follows from
  Theorem~\ref{thm:boglobalpar} that there is a path $\s_t \in
  \P(N,\ua)$ joining $\s_0,\s_1$ along which the bending angles are
  $\theta_{\a_i}(t) = t u_i + (1-t)v_i$ where $0 \leq t \leq 1$.  (The
  case where $u_i$ or $v_i$ is zero can be handled as in the proof of
  Proposition~\ref{prop:limit}.)  Note that $\theta_{\a_i}(t) \equiv
  u_i=v_i$ when $i>q$.

Along this path we have
$$\frac{d{l_{\a_i}}}{ dt}(\sigma_t)= \sum_{j=1}^n \frac{\partial
  l_{\a_i}}{\partial\theta_{\a_j}}(\sigma_t)
\,\frac{d\theta_{\a_j}}{dt}= \sum_{j=1}^n \frac{\partial
  l_{\a_i}}{\partial\theta_{\a_j}}(\sigma_t) (u_j-v_j).$$
If we
multiply both sides by $(u_i - v_i)$ and sum over $i$, it follows from
Corollary~\ref{cor:posdef} that $ \sum_{i=1}^n
(u_i-v_i)\frac{d{l_{\a_i}}}{ dt}(\sigma_t)
=  \sum_{i=1}^q
(u_i-v_i)\frac{d{l_{\a_i}}}{ dt}(\sigma_t) <0$.  Thus integrating
along $\sigma_t$ we find $\sum_{i=1}^q (u_i-v_i)(l_{\a_i}(\s_0) -
l_{\a_i}(\s_1) )<0$.  But this is impossible, since by our assumption
$l_{\a_i}(\s_0) = l_{\a_i}(\s_1)$ for all $i \leq q$.
\end{proof}

\noindent {\bf Theorem A (Length parameterization)}
{\em Let $L: \P(N,\ua) \rightarrow \RR^n$ be the map which associates
  to each structure $\sigma$ the hyperbolic lengths
  $(l_{\a_1}(\sigma), \ldots,l_{\a_n}(\s))$ of the curves in the
  bending locus $\ua=\{\a_1,\ldots,\a_n\}$. Then $L$ is an injective
  local diffeomorphism.}
\medskip

\noindent{\sc Proof.}
This is a special case of Theorem~\ref{cor:mixedparams}.
\ q.e.d. \vskip 0.2in

\begin{cor}
\label{cor:mixed}
If  $c_j \in [0,\pi]$ for $ j > q$, then
$$\{ \sigma \in \P(N,\ua)\,|\, \theta_{j}(\sigma)\equiv c_j,
j > q \}$$  is parameterized
by the lengths $l_{\a_j}$ with $ j \leq q$.
\end{cor}

We finish with a couple of other easy consequences of the
positive definiteness of $\mbox{d}L(\sigma)$.

\begin{cor} Suppose that $\sigma \in \P(N,\ua)$. Then for all
$\a_i \in \ua$,
$$\frac{\partial l_{\a_i}}{\partial \varphi_{\a_i}}(\sigma) > 0.$$
Therefore in the doubled cone manifold  $\Delta(\sigma)$, the length
of the singular locus is always increasing as a function
of the cone angle.
\end{cor}
\noindent{\bf Remark} For a general cone-manifold it is not true that
the derivative of the length of the singular locus with respect to
  the cone angle is always strictly positive.
  For an example in the case of the figure-eight
knot complement, see~\cite{choi}.

\begin{cor}
The volume of the convex core is a strictly concave function on
$\P(N,\ua)$
as a function of the bending angles,
with a global maximum at the unique structure for which all the bending
angles are $\pi$.
\end{cor}
\begin{proof} Let $\sigma_0 \in \P(N,\ua)$ and let $\sigma_t,\  t \in
[0,1]$ denote the
path of structures corresponding to the bending
angles $\theta_i(\sigma_t)= \theta_{i}(\sigma_0) + b_it$ where $b_i =
\pi-\theta_{i}(\sigma_0)$, as in Proposition~\ref{prop:limit}.
By Schl\"{a}fli's formula,
along this path we have
$$\frac{d{\rm Vol}}{ dt}(\sigma_t)=
   \sum_{i=1}^n l_{\a_i}(\sigma_t)\,\frac{d\theta_{i}}{dt}(\sigma_t)=
   \sum_{i=1}^n
b_i l_{\a_i}(\sigma_t)$$
which is strictly positive except at the unique maximal cusp
$\sigma_1$.

In a similar way, we can construct a linear path $\sigma_t$ between
$\sigma_0$ and any other point $ \sigma' \in \cal P(N,\ua)$ with
angles given by $\theta_i(\sigma_t)= \theta_{i}(\sigma_0) + c_it$,
where $c_i = \theta_i(\sigma') - \theta_{i}(\sigma_0)$.
    To prove concavity, we have to show the second derivative is
negative. We have:
$$ -\frac{d^2 {\rm Vol}}{ {d}t^2}(\sigma_t)
=
-\frac{d}{dt}\left(\sum_{i=1}^n c_i
l_{\a_i}(\sigma_t)\right)
=-\sum_{i,j=1}^n
\frac{\partial l_{\a_i}}{\partial\theta_{\a_j}}(\sigma_t)
c_ic_j = \frac{1}{2}\sum_{i,j=1}^n
\frac{\partial l_{\a_i}}{\partial\varphi_{\a_j}}(\sigma_t)
c_ic_j.
$$
The expression on the right is the value of the positive
definite quadratic form $\mbox{d}L(\sigma_t)$ evaluated on the vector
$(c_1,\ldots,c_d)$ and is therefore positive. The claim follows.
\end{proof}

\section*{Appendix}
\label{sec:appendix1}

Here is the Bonahon-Otal proof of
Proposition~\ref{prop:sufnonemptypleating}. The topological details
they omit are explained in Lemma~\ref{lemma:doubletopology}.\\

\noindent{\sc Proof of Proposition~\ref{prop:sufnonemptypleating}.}
Let $\bar M$ be the compact $3$-manifold obtained from the
double $D \bar N$ by removing disjoint tubular neighborhoods of
the curves in $\ua$ and let $M$ be its interior.
Lemma~\ref{lemma:doubletopology} shows that
the condition that
$\u{\a}$ be doubly incompressible with respect to $(\bar N,
\dd\bar N)$ is precisely the condition needed in order to apply
Thurston's hyperbolization theorem for Haken manifolds
(as in for example~\cite{Mor} p.52
or~\cite{Kap}  Theorem 1.42) to $\bar M$.
This gives a finite volume complete hyperbolic structure on $M$
in which every boundary component of $\bar M$ corresponds to a
rank-$2$ cusp.

There is a natural orientation reversing
involution $\tau:M \to M$ which interchanges the two copies
of $\bar N-\ua$. By Mostow rigidity, $\tau$ is homotopic to an
isometry, which we again denote by $\tau$. Taking the quotient of
$M$ by $\tau$ endows $\bar N - \ua$ with a
complete hyperbolic metric of finite volume. The boundary of
$\dd\bar N - \ua$ is the fixed point set of $\tau$ and therefore is
the union of totally geodesic surfaces. (The fixed point set of an
orientation reversing involution of $\HH^3$ is a plane.) Since all
curves in $\ua$ are parabolic by construction, we have exhibited a
convex structure on $(\bar N,\ua)$.
\ q.e.d. \vskip .2in

\renewcommand{\theappendthm}{A.\arabic{appendthm}}
\begin{appendlemma}
\label{lemma:doubletopology}
Let $\bar N$ be a compact orientable $3$-manifold whose interior
$N$ admits a complete hyperbolic structure. Assume that $\dd \bar N$
is non-empty and that it contains only surfaces of strictly negative
Euler characteristic Let
$\u{\a} = \{\a_1, \ldots, \a_{n}\}$ be a doubly incompressible curve
system on $\dd\bar N$ and let $M$ be as defined above. Then
$\dd \bar M$ is incompressible in $\bar M$ and $\bar M$ is irreducible
and atoroidal.
\end{appendlemma}
     \begin{proof}
We begin with some notation. For each $i$,  let $U_{\a_i} $ be a regular
neighborhood of $\a_i$ in $\dd \bar N$ and let $\dd_0 \bar N = \dd N
- \cup_i U_{\a_i}$.
     Let $  N^+$ be the manifold with boundary obtained by removing, for
each $i$,  an open regular tubular neighborhood of $\a_i$
from $\bar N$ which intersects $\dd \bar N$ in $U_{\a_i}$.
Let $  N^- = \sigma(N^+)$. Then $N^+$ and $N^-$ are both homeomorphic
to $\bar N$. Clearly we can take $\bar M$ to be the union
of the images of $ N^+ $ and $ N^-$ in $ M$ glued along $\dd_0 \bar N$.
We refer to $N^+$ and $N^-$ as the {\em sides} of $M$. The assumption
that $N$ is hyperbolic implies that $N^\pm$ are irreducible and
atoroidal.

First we check that $\bar M$ is irreducible. Let $S$ be an embedded
sphere in $\bar M$. We may assume that $S \cap \dd \bar M =
\emptyset$, since the existence of a collar neighborhood of $\dd \bar
M$ in $\bar M$ allows $S$ to be pushed off $\dd \bar M$ if
necessary. Since the manifolds $N^\pm$ are irreducible,
we need only consider $S$ which intersects $\dd_0 \bar N$.
Assume that $S$ has been homotoped so that the intersection is
transverse. Then $S\, \cap\, \dd_0 \bar N$ is a finite union of
disjoint circles. Choose an innermost circle $C$, meaning that $C$ is
the boundary of a disk $D(C)$ in $S$ such that $\mbox{Int}\,D(C)\, \cap
\,\dd_0 \bar N = \emptyset$. Then $D(C)$ is completely
contained in one side of $\bar M$. Since $C \,\cap \,\ua = \emptyset $,
the condition ({\em D.2}\,),  that
every essential disk in $\bar N$ intersects $\ua$ at least $3$
times, implies that $D(C)$ is not essential. Therefore, $D(C)$
can be homotoped to a disk $D'(C)$ in $\dd_0 \bar N$ by a homotopy
fixing $C$ and then
pushed off $\dd_0 \bar N$ so that the number of circles in
$S \cap \dd_0 \bar N$ is reduced by one. By successively
applying the above process, we can homotope $S$ so that it no
longer intersects $\dd_0 \bar N$, which implies that it is
contained in one side of $\bar M$ and therefore bounds a $3$-ball.

Note that exactly the same arguments show that   $\dd_0 \bar N$ is
incompressible in $\bar M$.

Now we check that $\partial \overline{M}$ is incompressible. Recall
that the only irreducible orientable $3$-manifold with a
compressible torus boundary component is the solid torus. For suppose
that
a torus boundary component $T$ of an irreducible $3$-manifold $W$
contains
a non-trivial loop $C$ that bounds a disk $D$ in $W$. Let $D\times I$ be
a
product neighborhood of $D$ meeting $\partial W$ in $\partial D\times
I$. The $2$-sphere $D\times \{0,1\}\cup \overline{\partial T-\partial
D\times I}$ must bound a $3$-ball $B$ in $W$. The $3$-ball cannot
contain
$D\times I$ since then it would meet $\partial W$ in more than
$\overline{T-\partial D\times I}$. Therefore $W$ is the solid torus
obtained by attaching $D\times I$ to $B$ along $D\times \{0,1\}$. So, if
any component of $\partial \overline{M}$ were compressible,
$\overline{M}$ would be a solid
torus. Since $\partial_0\overline{N}$ is incompressible, it would
consist
of annuli,
contradicting the fact  that each of its components has negative
Euler characteristic.

Lastly, we check that $\bar M$ is atoroidal. Suppose there is
a $\ZZ \oplus \ZZ$ subgroup of $\pi_1(\bar M)$ which is
non-peripheral. The torus theorem
states that either there is an embedded incompressible torus
which is not boundary parallel or $\bar M$ is a Seifert fibered space,
see~\cite{BonG} Theorem 3.4. (Since $M$ is orientable, we do not need to
consider Klein bottles.)
We claim $\bar M$ cannot be Seifert fibered. As noted
   above,   $\dd_0 \bar N$ is incompressible in $\bar
M$. If $\bar M$ were Seifert fibered, then $\dd_0 \bar N$ would be
isotopic to a vertical or horizontal surface, see Theorem VI.34
~\cite{jaco}. In the first case,
$N^+$ must admit a Seifert fibering, but then $\dd N^+$ and hence
$\dd \bar N$ would consist of tori. In the second case, each
component of $\bar M$ cut along $\dd_0 \bar N$ is an $I$-bundle,
and the original $\bar N$ must have contained an annulus which
violates condition ({\em D.1}\,) or $\bar N$ is an $I$-bundle over a
pair of pants which
violates condition ({\em D.2}\,).

Now suppose there is an
incompressible torus $T$ embedded in $\bar M$. We may assume that $T
\cap \dd
\bar M =\emptyset$.
Since the manifolds $N^\pm$ are atoroidal, we may assume that
$T$ intersects $\dd_0 \bar N$ and that the intersection is transverse.
Then $T \cap \dd_0 \bar N$ is a disjoint union of circles. By the
same process as above, we can homotope $T$ to eliminate any circles
which are trivial in $T$. Since $T$ is a torus, the remaining
circles of intersection must all be parallel and
non-trivial. Furthermore, since $T$ is incompressible, the circles
are also non-trivial in $\dd_0 \bar N$.
    Take two adjacent
circles $C, C'$. Then the annulus $A$ they bound in $T$ is completely
contained in one side of $\bar M$. The condition that there
are no essential annuli in $\bar N$ with boundary in $\dd \bar N -
\ua$ implies that $A$ is not essential. Therefore $A$ can be homotoped
into an annulus $A'$ in $\dd N^+$ or $\dd N^-$ by a homotopy fixing
$C,C'$. If $A' \cap \dd(\dd_0 \bar N) = \emptyset$, then $A'$ is
contained in a component of $\dd_0 \bar N$. In this case, $A'$ can be
pushed off $\dd_0 \bar N$, thereby removing the circles of intersection
$C$ and $C'$. In this way, we can homotope $T$ so that no pair
of adjacent circles in $\dd_0 \bar N \cap T$ bound an annulus which
can be homotoped into a component of $\dd_0 \bar N$. Note that we
cannot eliminate all the circles in $\dd_0 \bar N \cap T$, since
this would imply that there is an incompressible torus entirely
contained in one side of $\bar M$. Also note that $\dd_0 \bar N
\cap T$ cannot contain only one circle, since $\dd N^\pm$ separates
$\bar M$ into two components. Now if $C,C'$ are two adjacent
circles in $\dd_0 \bar N \cap T$, they bound an annulus $A$
homotopic to an annulus $A'$ in $\dd N^+$ such that $A' \cap \dd(\dd_0
\bar N) \neq \emptyset$. Since no two curves in $\ua$ are freely
homotopic to one another, $A' \cap \dd(\dd_0 \bar N)$
consists of the two boundary curves of $U_{\a_i}$ for some $\a_i \in
\ua$. Since this is true of every annulus in $T$, it must be that
$T$ is homotopic into $\dd \bar M$. A  short further analysis
confirms that since $T$ is embedded, in fact there must be exactly
two distinct
circles $C,C'$ and hence that  $T$ is actually parallel to $\dd \bar M$.
\end{proof}
\small{
\end{document}